\documentclass[pdflatex,sn-mathphys-num]{sn-jnl}

\usepackage{graphicx}%
\usepackage{multirow}%
\usepackage{amsmath,amssymb,amsfonts}%
\usepackage{amsthm}%
\usepackage{mathrsfs}%
\usepackage[title]{appendix}%
\usepackage{xcolor}%
\usepackage[table]{xcolor}
\usepackage{textcomp}%
\usepackage{manyfoot}%
\usepackage{booktabs}%
\usepackage{algorithm}%
\usepackage{algorithmicx}%
\usepackage{booktabs}
\usepackage{algpseudocode}%
\usepackage{listings}%
\usepackage{bm}
\usepackage{cleveref}
\usepackage{subcaption}
\usepackage{placeins}
\usepackage{tikz}
\usepackage{pgf}
\usepackage{adjustbox}
\usetikzlibrary{arrows.meta, decorations.pathreplacing, patterns}
\crefname{subsection}{subsection}{subsections}
\Crefname{subsection}{Subsection}{Subsections}
\usepackage{comment}
\usepackage[style=authoryear,sorting=nyt,giveninits=true,backend=biber]{biblatex}
\let\cite\parencite

\numberwithin{equation}{section}
\newcommand{\GhatLIS}[0]{\mathbf{\widehat{G}}_{\textLISMR}}
\newcommand{\GSpa}[0]{\mathbf{G}_{\text{\tiny OLR}}}
\newcommand{\proBasis}[0]{\mathbf{W}}
\newcommand{\approxBasis}[0]{\mathbf{V}}
\newcommand{\textLISMR}[0]{\text{\tiny LIS}}

\definecolor{lightgray}{gray}{0.9}
\usepackage{nicematrix}
\usepackage{comment}

\theoremstyle{thmstyleone}%
\theoremstyle{thmstyletwo}%

\theoremstyle{thmstylethree}%

\begin{document}

\title[Article Title]{Likelihood-informed Model Reduction for Bayesian Inference of Static Structural Loads}


\author*[1,4]{\fnm{Jakob} \sur{Scheffels}}\email{jakob.scheffels@tum.de}

\author[2,4]{\fnm{Elizabeth} \sur{Qian}}\email{eqian@gatech.edu}

\author[1]{\fnm{Iason} \sur{Papaioannou}}\email{iason.papaioannou@tum.de}

\author[3]{\fnm{Elisabeth} \sur{Ullmann}}\email{elisabeth.ullmann@ma.tum.de}

\affil[1]{\orgdiv{Engineering Risk Analysis Group}, \orgname{TUM School of Engineering and Design, Technical University of Munich}, \orgaddress{\street{Arcisstr 21}, \city{Munich}, \postcode{80333}, \state{Bavaria}, \country{Germany}}}

\affil[2]{\orgdiv{School of Aerospace Engineering and School of Computational Science and Engineering}, \orgname{Georgia Institute of Technology}, \orgaddress{\street{North Avenue}, \city{Atlanta}, \postcode{30332}, \state{GA}, \country{USA}}}

\affil[3]{\orgdiv{Department of Mathematics}, \orgname{TUM School of Computation, Information and Technology, Technical University of Munich}, \orgaddress{\street{Boltzmannstr 3}, \city{Garching}, \postcode{85748}, \state{Bavaria}, \country{Germany}}}

\affil[4]{\orgdiv{Institute for Advanced Study}, \orgname{Technical University of Munich}, \orgaddress{\street{Lichtenbergstr 2A}, \city{Garching}, \postcode{85748}, \state{Bavaria}, \country{Germany}}}

\abstract{
Bayesian inverse problems use data to update a prior probability distribution on uncertain parameter values to a posterior distribution. Such problems arise in many structural engineering applications, but computational solution of Bayesian inverse problems is often expensive because standard solution approaches require many evaluations of the forward model mapping the parameter value to predicted observations. In many settings, this forward model is expensive because it requires the solution of a high-dimensional discretization of a partial differential equation. However, Bayesian inverse problems often exhibit low-dimensional structure because the available data are primarily informative (relative to the prior) in a low-dimensional subspace, sometimes called the likelihood-informed subspace (LIS). 
This paper proposes a new projection-based model reduction method for static linear systems that exploits this low-dimensional structure in the setting where the unknown parameter is the right-hand-side forcing\textcolor{black}{, giving rise to a linear inverse problem}. 
The proposed method projects the governing partial differential equation onto the likelihood-informed subspace, yielding a computationally efficient reduced model that can be used to accelerate the solution of the inverse problem \textcolor{black}{and subsequent downstream computations}. Numerical experiments on two structural engineering model problems demonstrate that the proposed approach can successfully exploit the intrinsic low-dimensionality of the problem, obtaining relative errors \textcolor{black}{in} $\mathcal{O}(10^{-10})$ in the inverse problem solution with a 10$\times$ \textcolor{black}{and} 100$\times$ lower-dimensional model\textcolor{black}{, respectively}.
}

\keywords{Bayesian inference, model order reduction, likelihood-informed subspace, posterior approximation, structural engineering}



\maketitle

\section{Introduction}\label{sec1}
Uncertain parameters in structural engineering models may arise due to a variety of factors, including heterogeneity of materials, manufacturing variations, and environmental influences. 
The performance of a structure can be highly sensitive to these uncertainties in the model parameters, making estimation of these parameters crucial for structural design \cite{hurtadoBayesian5StoryBuilding} and digital twin modeling \cite{arconesDigitalTwinBridges,torzoniDigitalTwinFramework,thelenComprehensiveReviewDigitalTwin} to assess the ongoing health of a structure. For example, in health monitoring of a subway tunnel segment, uncertain model parameters include the tunnel material properties, ground support, and live load acting on the tunnel \cite{huangTunnel}. These parameters must be accurately estimated  for accurate modeling of the structure. 

Inverse problems describe the estimation of an unknown parameter from observed data. 
The Bayesian approach~\cite{stuartInverseProblemsBayesian2010b} to inverse problems endows the uncertain parameters with a prior probability distribution that encodes the baseline uncertainty in the parameter value, and then uses observed data to update the prior to a posterior distribution.
However, many model parameters, such as the Young's modulus or the load acting on the structure, cannot be measured directly. Instead, we are limited to indirect observations of their effect, such as the resulting displacement of the structure, rather than the parameter itself. For such problems, solving inverse problems requires repeated evaluation of the \textit{forward model} that maps the uncertain parameters to the measurable quantities~\cite{auvinenLARGESCALEKALMANFILTERING, martinStochasticNewtonMCMC2012a}. In situations where the forward model is expensive, e.g.\ when it entails the solution of a high-dimensional discretization of a partial differential equation (PDE), computational solution of inverse problems can be prohibitively expensive. 

To reduce the computational cost of inference, surrogate modeling methods aim to cheaply approximate the forward model. One approach to building surrogates for Bayesian inverse problems is polynomial chaos expansion (PCE)~\cite{AdaptiveMultifidelityPCE2019,novakPhysicsinformedPolynomialChaos2024}, which approximates the forward model in the span of a truncated set of cheap-to-evaluate basis functions. However, the truncation of the basis can lead to high errors due to inexpressivity~\cite{LimitationsPolynomialChaos2015}.
In contrast, universal function approximators such as Gaussian processes and neural networks~\cite{dinkelSolvingBayesianInverse2023, villaniAdaptiveGaussianProcess2024, yanAdaptiveSurrogateModeling2020, deveneyDeepSurrogateApproach2021,pasparakisBayesianNeuralNetworks2025, pfortnerPhysicsInformedGaussianProcess2024} are popular surrogate modeling approaches due to their expressivity. These methods are highly flexible and can accelerate inference computations by orders of magnitude, but often require large volumes of training data, which may not be available in practical applications. 

Projection-based model reduction describes a class of surrogate modeling approaches which obtain cheap reduced models by projecting the governing PDE onto a low-dimensional subspace \cite{Antoulas2005, benner2015survey}. The most common projection-based model reduction method used in the context of Bayesian updating is proper orthogonal decomposition (POD) \cite{lumley1981coherent,sirovich1987turbulence, ghattasLearningPhysicsbasedModels2021a, nguyenModelOrderReduction2014, cuiDatadrivenModelReduction2015, xiongAcceleratingBayesianInference2021, raoInverseParameterEstimation2024}. To obtain the reduced-order model, POD approximates the system state in the span of the leading principal components of available state snapshot data; these data may either be given or obtained by solving the high-dimensional model for samples of the uncertain parameters. The reduced basis method~\cite{hesthaven2022reduced,rozza2024short,chenSteinVariationalReduced2021,silvaReducedBasisEnsemble2023} extends POD to parametrized problems. 
While the snapshots used to define the projection basis are often chosen independently of the Bayesian inverse problem context~\cite{nguyenModelOrderReduction2014}, POD and RB methods may be tailored to Bayesian inverse problems by sampling snapshot data from either the prior \cite{ghattasLearningPhysicsbasedModels2021a} or posterior distribution \cite{cuiDatadrivenModelReduction2015}.
However, the primary goal of these snapshot-based methods is accurate reconstruction of the available snapshot data, which may be expensive to generate. Additionally, accurate reconstruction of the full state may not be necessary for accurate approximation of inverse problem solutions, because observation data are often primarily informative (relative to the prior) only in a low-dimensional subspace of the full state space~\cite{spantiniOptimalLowRank}.

This low-dimensional structure of inverse problems is well-known and has been exploited by a variety of \textit{dimension reduction} approaches that achieve computational savings by only approximating the posterior distribution within the low-dimensional space that is informed by the data, sometimes called the likelihood-informed subspace (LIS)~\cite{spantiniOptimalLowRank, bui-thanhExtremescaleUQBayesian2012, bui-thanhComputationalFrameworkInfiniteDimensional2013,cuiLikelihoodinformedDimensionReduction2014b,zahm2022certified}.
For linear Gaussian Bayesian inverse problems, \cite{spantiniOptimalLowRank} shows that such approximations are optimal. 
However, these dimension reduction approaches generally still require the evaluation of the high-dimensional forward model for parameter values within the LIS. \textcolor{black}{Additionally, any downstream computations requiring evaluation of the model (e.g., uncertainty propagation of the parameter posterior to the state space) still must be done using the full model in the high-dimensional space}. To achieve additional computational acceleration, several recent works have therefore developed new projection-based model reduction approaches that exploit the low-dimensional LIS for Bayesian inverse problems to obtain efficient low-dimensional reduced models~ \cite{qianModelReductionLinearBalancing, josieLowRankPrior, konigTimeLimitedBalancedTruncation2023,stavrinidesEnsembleKalmanApproach2025, freitagInferenceOrientedBalancedTruncation2024a}. However, these approaches focus on model reduction for inferring the initial state of linear dynamical systems.

In this work, we propose a new LIS-based model reduction approach applicable to linear \textit{static} systems where the parameter to be inferred is the unknown right-hand side of the system, and the observed data come from noisy linear measurements of the state or system response. \textcolor{black}{Our approach defines reduced representations of both the parameter and the system state, enabling both efficient approximate solution of the inverse problem within the reduced space and also efficient downstream uncertainty propagation computations that can be performed with the reduced model instead of the high-dimensional full model.} We provide two numerical experiments on structural engineering model problems demonstrating that the proposed method leads to more accurate results than state-of-the-art model reduction methods. 
\Cref{sec:background} introduces the inverse problem we consider and background on model and dimension reduction approaches for inverse problems. \Cref{sec:LIS} introduces our new method \textcolor{black}{and compares cost and storage requirements of the proposed method with existing approaches}. In \Cref{sec:Examples}, we apply the method to inferring the unknown load on a cantilever bar as well as a tunnel segment. \Cref{sec:Conclusion} 
\textcolor{black}{contains concluding remarks}.

\section{Background}\label{sec:background}
\textcolor{black}{We introduce simulation-based inverse problems in \Cref{sec:simIP}, together with a discussion of the computational challenges in their solution that motivate this work. \Cref{sec:linGauIP} then introduces the specific static structural load inference problem that we consider in this work.} \Cref{sec:MOR} introduces projection-based model reduction and the POD method. \Cref{sec:Spantini} introduces optimal low-rank approximations of linear Gaussian Bayesian inverse problems from~\cite{spantiniOptimalLowRank}.

\subsection{Simulation-based inverse problems}\label{sec:simIP}
Consider a computational model
\begin{equation}
    \mathbf{u}=\mathbf{M}(\mathbf{f})\label{eq:compModel},
\end{equation}
where $\mathbf{u}\in\mathbb{R}^d$ is the system state, $\mathbf{M}:\mathbb{R}^p\rightarrow\mathbb{R}^d$ is a nonlinear model mapping the parameters to the system state and $\mathbf{f}\in\mathbb{R}^p$ is an unknown parameter that must be estimated or inferred from measurements $\mathbf{y}\in\mathbb{R}^m$ given by 
\begin{equation}
    \mathbf{y} = \mathbf{C}(\mathbf{u})+\bm\epsilon, \label{eq:generalIP}
\end{equation}
where $\mathbf{C}:\mathbb{R}^d\rightarrow\mathbb{R}^m$ is a map from the $d$-dimensional state to the $m$-dimensional measurements and $\bm\epsilon\in\mathbb{R}^{m}$ is the measurement noise. 
We are interested in the setting where the model $\mathbf{M}(\mathbf{f})$ involves the solution of an underlying PDE that characterizes the system's response to the unknown parameters $\mathbf{f}$. 

The inverse problem is then fully defined by the observation model,
\begin{align*}
    \mathbf{y} = \mathbf{G}(\mathbf{f}) + \bm\epsilon,
\end{align*}
where the forward map $\mathbf{G}$ relating the parameter $\mathbf{f}$ to measurements $\mathbf{y}$ is given by the composition of the PDE-based model and the measurement operator:
\begin{equation}
    \mathbf{G}=\mathbf{C}\circ\mathbf{M}\label{eq:generalG}.
\end{equation}
Note that predicting measurements for given parameter values requires evaluating the computational model in \eqref{eq:compModel}, which is often expensive.

We use a Bayesian approach to solve the desired inverse problem. The approach treats the measurements $\mathbf{y}$, the parameters $\mathbf{f}$ and the measurement noise $\bm\epsilon$ as random variables. We assign a prior distribution $\pi_\mathbf{F}(\mathbf{f})$ to the unknown parameters   that describes the uncertainty of $\mathbf{f}$ prior to obtaining the measurements.
The conditional distribution of $\mathbf{f}$ given the measurements $\mathbf{y}$, or  posterior distribution of $\mathbf{f}$, is defined following Bayes' rule as being proportional to the product of the likelihood of the observations $\pi_{\mathbf{Y}|\mathbf{F}}(\mathbf{y}|\mathbf{f})$ and the prior distribution
\begin{equation*}
    \pi_{\mathbf{F}|\mathbf{Y}}(\mathbf{f}|\mathbf{y})\propto \pi_{\mathbf{Y}|\mathbf{F}}(\mathbf{y}|\mathbf{f})\pi_\mathbf{F}(\mathbf{f}).
\end{equation*}
In general, there does not exist an analytical solution to obtain $\pi_{\mathbf{F}|\mathbf{Y}}(\mathbf{f}|\mathbf{y})$ and we have to use sampling methods to sample from the posterior distribution \cite{beck2002bayesian, martinStochasticNewtonMCMC2012a}. These methods require multiple evaluations of the likelihood of the data. Because each likelihood evaluation involves the solution of a high-dimensional computational model, these methods are computationally expensive.
\textcolor{black}{Reduced models, such as those introduced in \Cref{sec:MOR} and proposed in this work, have the potential to make characterizing Bayesian posteriors more tractable. In order to lay a foundation for inference-oriented reduced modeling methods, this work develops our method for linear inverse problems, as specifically described in the next section.}

\subsection{\textcolor{black}{Bayesian inference of static structural loads}}\label{sec:linGauIP}
In this paper, we focus on the 
\textcolor{black}{inference of static structural loads}, where both the observation operator and the PDE model are linear \textcolor{black}{and uncertainties only lie on the right-hand-side of the PDE}:
\begin{equation}
    \mathbf{y}=\mathbf{C}\mathbf{u}+\bm\epsilon,\quad \text{s.t.:} \quad \mathbf{K}\mathbf{u}=\mathbf{f}\label{eq:linearPDE},
\end{equation}
where $\mathbf{K}\in\mathbb{R}^{d\times d}$ is a linear system matrix and thus $\mathbf{M}(\mathbf{f}) = \mathbf{K}^{-1}\mathbf{f}$. For example, such a system arises in the inference of unknown static structural loads, where $\mathbf{u}$ is the displacement, $\mathbf{K}$ the stiffness matrix and $\mathbf{f}$ the unknown load. The corresponding inverse problem is given by
\begin{flalign}
    \mathbf{y}=\mathbf{G}\cdot \mathbf{f}+\bm\epsilon \label{eq:linearIP},\\
    \mathbf{G} = \mathbf{C}\cdot \mathbf{K}^{-1}, \label{eq:forwardOperator}
\end{flalign}
where $\mathbf{G}\in\mathbb{R}^{m\times d}$ is now a linear forward operator.

We assume a Gaussian prior distribution $\pi_\mathbf{F}(\mathbf{f})=\mathcal{N}(\bm\mu,\bm\Gamma)$ with known mean and covariance $\bm\mu\in\mathbb{R}^d$ and $\mathbf{\Gamma}\in\mathbb{R}^{d\times d}$, and Gaussian measurement noise $\bm\epsilon\sim\mathcal{N}(\mathbf{0},\bm\Gamma_{\text{obs}})$ with $\mathbf{\Gamma}_{\text{obs}}\in\mathbb{R}^{m\times m}$. Then, the posterior distribution $\pi_{\mathbf{F}|\mathbf{Y}}(\mathbf{f}|\mathbf{y})$ will also be Gaussian and is fully described by the posterior mean $\bm{\mu}_{\text{pos}}\in\mathbb{R}^d$ and covariance $\bm\Gamma_{\text{pos}}\in\mathbb{R}^{d\times d}$ whose analytical expressions are given by~\cite{stuartInverseProblemsBayesian2010b} 
\begin{flalign}
    \bm\mu_{\text{pos}}&=\bm\mu+\bm\Gamma \mathbf{G}^\top(\mathbf{G}\mathbf{\Gamma}\mathbf{G}^\top+\mathbf{\Gamma}_{\text{obs}})^{-1}(\mathbf{y}-\mathbf{G}\bm\mu)\label{eq:posMean},\\
    \bm\Gamma_{\text{pos}}&=\bm\Gamma-\bm\Gamma \mathbf{G}^\top(\mathbf{G}\bm\Gamma 
    \mathbf{G}^\top+\bm\Gamma_{\text{obs}})^{-1}\mathbf{G}\bm\Gamma\label{eq:posCov}.
\end{flalign}
Note that computing the posterior statistics requires evaluation of the high-dimensional forward model $\mathbf{G}$ and its adjoint, which can be expensive. This challenge is exacerbated in settings where the solution to the inverse problem must be computed repeatedly for different data, e.g., in digital twin settings. To reduce the cost of inference in these settings, model reduction methods can be employed.

\subsection{Projection-based model reduction}\label{sec:MOR}
Projection-based model reduction builds a reduced-order model by Petrov-Galerkin projection of the PDE into a low-dimensional subspace. Let $\approxBasis\in\mathbb{R}^{d\times r},\proBasis\in\mathbb{R}^{d\times r}$ denote a trial and test basis, respectively. Projection-based model reduction methods approximate the state in the span of the trial basis 
\begin{equation}
    \mathbf{u}\approx\mathbf{V}\cdot \mathbf{\widehat{u}},\label{eq:approximationFullState}
\end{equation}
where $\mathbf{\widehat{u}}\in\mathbb{R}^r$ is the reduced state. Inserting \eqref{eq:approximationFullState} into the discretized PDE given in \eqref{eq:linearPDE} and enforcing the Petrov-Galerkin orthogonality condition that the residual is orthogonal to the test basis leads to
\begin{flalign*}
    \mathbf{W}^\top \mathbf{K}(\mathbf{V}\mathbf{\widehat{u}})=\mathbf{W}^\top \mathbf{f}.
\end{flalign*}
Defining $\mathbf{\widehat{K}}=\mathbf{W}^\top \mathbf{K}\mathbf{V}\in\mathbb{R}^{r\times r}$ and $\mathbf{\widehat{f}}=\mathbf{W}^\top \mathbf{f}\in\mathbb{R}^r$, we obtain the reduced model
\begin{equation}
    \mathbf{\widehat{K}}\mathbf{\widehat{u}}=\mathbf{\widehat{f}},\label{eq:reducedPDE}
\end{equation}
where solving for $\mathbf{\widehat{u}}$ only requires the inversion of $\mathbf{\widehat{K}}\in\mathbb{R}^{r\times r}$ instead of $\mathbf{K}\in\mathbb{R}^{d\times d}$. 

In the context of linear inverse problems, we can use the reduced model in \eqref{eq:reducedPDE} to formulate an approximation of \eqref{eq:linearIP} defined by 
\begin{equation}
    \mathbf{y}=\mathbf{\widehat{C}}\mathbf{\widehat{u}}+\bm\epsilon,\quad \text{s.t.:} \quad \mathbf{\widehat{K}}\mathbf{\widehat{u}}=\mathbf{\widehat{f}},\label{eq:reducedIP}
\end{equation}
where $\mathbf{\widehat{C}}\in\mathbb{R}^{m\times r}$ is the reduced observation operator obtained by $\mathbf{\widehat{C}}=\mathbf{C}\approxBasis$.  
The $\widehat{\phantom{x}}$ notation denotes a reduced, low-dimensional operator. Typically the high-dimensional computations that define $\approxBasis$ and $\proBasis$ \textcolor{black}{as well as computations of} the reduced quantities $\widehat{\mathbf{C}}$ 
\textcolor{black}{and} $\widehat{\mathbf{K}}$ 
are viewed as occurring during an ``offline'' phase in which the reduced model is obtained. \textcolor{black}{As new data are obtained in the ``online" phase, the reduced model \eqref{eq:reducedPDE} can be used to cheaply approximate the inverse problem.} 
 A key feature of model reduction methods is that the computational cost of the online phase scales only with the reduced dimension $r$, and is independent of the original high dimension $d$.

For the general case of a Petrov-Galerkin projection, the test basis $\proBasis$ and the trial basis $\approxBasis$ differ. If $\proBasis=\approxBasis$, we speak of a Galerkin projection. By far the most commonly used model reduction method is Galerkin projection using a proper orthogonal decomposition (POD) basis. 
POD-based model order reduction approximates the high-dimensional computational model based on a set of snapshots of the full state. For the specific problems we consider, the snapshots arise from sampling $\pi_\mathbf{F}(\mathbf{f})$ and evaluating the PDE for these samples to obtain the corresponding displacements. Assuming that we have generated $N$ samples of the random parameters $[\mathbf{f}_1,...,\mathbf{f}_N]$, we collect the resulting displacements given by evaluating \eqref{eq:linearPDE} for each snapshot in a matrix $\mathbf{U}=[\mathbf{u}_1,...,\mathbf{u}_N]\in\mathbb{R}^{d\times N}$. Let $\mathbf{U}=\bm\Phi\bm\Sigma \bm\Psi^\top$ denote the singular value decomposition of the snapshot matrix. Then, the leading $r$ left singular vectors contained in the submatrix $\bm\Phi_r\in\mathbb{R}^{d\times r}$ are a basis for the optimal rank-$r$ linear subspace for reconstructing the snapshot data~\cite{eckartApproximationOneMatrix1936}. Galerkin-POD models are obtained by setting $\mathbf{W}=\mathbf{V}=\bm\Phi_r$. The approximated inverse problem \eqref{eq:reducedIP} using a POD-reduced model is given by
\begin{flalign}
    \mathbf{y}&=\mathbf{\widehat{C}}_{\text{\tiny POD}}\mathbf{\widehat{u}}_{\text{\tiny POD}}+\bm\epsilon,\quad \text{s.t.:} \quad \mathbf{\widehat{K}}_{\text{\tiny POD}}\mathbf{\widehat{u}}_{\text{\tiny POD}}=\mathbf{\widehat{f}}_{\text{\tiny POD}},\nonumber\\
    \mathbf{y}&=\mathbf{\widehat{G}}_{\text{\tiny POD}}\mathbf{\widehat{f}}_{\text{\tiny POD}}+\bm\epsilon, \label{eq:PODIP}
\end{flalign}
where $\mathbf{\widehat{G}}_{\text{\tiny POD}}=\mathbf{C}\bm{\Phi}_r\mathbf{\widehat{K}}_{\text{\tiny POD}}^{-1}\in\mathbb{R}^{m\times r}$ is the approximation of the forward operator using the POD-reduced model $\mathbf{\widehat{K}}_{\text{\tiny POD}}=\bm\Phi_r^\top\mathbf{K}\bm\Phi_r$ and $\mathbf{\widehat{f}}_{\text{\tiny POD}}=\bm{\Phi}_r^\top\mathbf{f}\in\mathbb{R}^r$ is the reduced parameter. The reduced parameter $\mathbf{\widehat{f}}_{\text{\tiny POD}}$ has Gaussian prior distribution with $\bm{\widehat{\mu}}_{\text{\tiny POD}}=\bm\Phi_r^\top\bm\mu\in\mathbb{R}^{r}$ and $\mathbf{\widehat{\Gamma}}_{\text{\tiny POD}}=\bm\Phi_r^\top\mathbf{\Gamma}\bm\Phi_r\in\mathbb{R}^{r\times r}$ leading to a Gaussian posterior defined by
\begin{flalign}
    \bm{\widehat{\mu}}_{\text{pos}}^{\text{\tiny POD}}&=\bm{\widehat{\mu}}_{\text{\tiny POD}}+\mathbf{\widehat{\Gamma}}_{\text{\tiny POD}} \mathbf{\widehat{G}}_{\text{\tiny POD}}^\top(\mathbf{\widehat{G}}_{\text{\tiny POD}}\mathbf{\widehat{\Gamma}}_{\text{\tiny POD}} \mathbf{\widehat{G}}_{\text{\tiny POD}}^\top+\mathbf{\Gamma}_{\text{obs}})^{-1}(\mathbf{y}-\mathbf{\widehat{G}}_{\text{\tiny POD}}\bm{\widehat{\mu}}_{\text{\tiny POD}})\label{eq:posMeanPODReduced},\\
    \mathbf{\widehat{\Gamma}}_{\text{pos}}^{{\text{\tiny POD}}}&=\mathbf{\widehat{\Gamma}}_{\text{\tiny POD}}-\mathbf{\widehat{\Gamma}}_{\text{\tiny POD}} \mathbf{\widehat{G}}_{\text{\tiny POD}}^\top(\mathbf{\widehat{G}}_{\text{\tiny POD}}\mathbf{\widehat{\Gamma}}_{\text{\tiny POD}} \mathbf{\widehat{G}}_{\text{\tiny POD}}^\top+\bm\Gamma_{\text{obs}})^{-1}\mathbf{\widehat{G}}_{\text{\tiny POD}}\mathbf{\widehat{\Gamma}}_{\text{\tiny POD}},\label{eq:posCovPODReduced}
\end{flalign}
where $\bm{\widehat{\mu}}_{\text{pos}}^{\text{\tiny POD}}\in\mathbb{R}^r$ and $\mathbf{\widehat{\Gamma}}_{\text{pos}}^{\text{\tiny POD}}\in\mathbb{R}^{r\times r}$ are reduced representations of the posterior statistics. The corresponding high-dimensional approximation of the posterior mean of $\mathbf{f}$ is given by
\begin{equation}
    \bm\mu_{\text{pos}}^{\text{\tiny POD}}=\textcolor{black}{(\mathbf{I}_d-\bm\Phi_r\bm\Phi_r^\top)\bm\mu+}\bm\Phi_r\bm{\widehat{\mu}}_{\text{pos}}^{\text{\tiny POD}}\in\mathbb{R}^d,
\end{equation}
and the low-rank update posterior covariance approximation is given by
\begin{equation}
    \mathbf{\Gamma}_{\text{pos}}^{\text{\tiny POD}}=\mathbf{\Gamma}-\mathbf{\Gamma}\bm\Phi_r\mathbf{\widehat{G}}_{\text{\tiny POD}}^\top(\mathbf{\widehat{G}}_{\text{\tiny POD}}\bm\Phi_r^\top\mathbf{\Gamma}\bm\Phi_r\mathbf{\widehat{G}}_{\text{\tiny POD}}^\top+\mathbf{\Gamma}_{\text{obs}})^{-1}\mathbf{\widehat{G}}_{\text{\tiny POD}}^\top\bm\Phi_r^\top\mathbf{\Gamma}\in\mathbb{R}^{d\times d},\label{eq:posCovPODMapped}
\end{equation}
where the prior is kept in the subspace orthogonal of the POD basis. Using \eqref{eq:posMeanPODReduced}--\eqref{eq:posCovPODMapped}, we can solve an $r$-dimensional inverse problem to obtain an approximation of the $d$-dimensional posterior quantities. The reduced forward operator $\mathbf{\widehat{G}}_{\text{\tiny POD}}$ has storage and evaluation costs \textcolor{black}{in} $\mathcal{O}(mr)$ instead of $\mathcal{O}(md+d^2)$ leading to significant computational savings for $r\ll d$. \textcolor{black}{Note that the full forward operator has lower storage costs if $\mathbf{K}$ is sparse banded.} 

For POD-reduced models, the bases $\mathbf{W}=\mathbf{V}$ depend on the samples used to generate $\mathbf{U}$. In the context of Bayesian inverse problems, snapshots can either be generated by sampling from the prior \cite{ghattasLearningPhysicsbasedModels2021a} or adaptively from the posterior distribution \cite{cuiDatadrivenModelReduction2015}, or are based on given data \cite{nguyenModelOrderReduction2014}. However, POD does not explicitly consider the transition from prior to posterior distribution when computing $\bm\Phi_r$, but focuses on the reconstruction of the full state, which might not be necessary for an accurate approximation of the solution of an inverse problem. Additionally, the data-driven basis calculation requires $N$ evaluations of the full computational model for generating snapshots first, leading to a large computational overhead.

\subsection{Optimal low-rank approximation for linear Gaussian Bayesian inverse problems}\label{sec:Spantini}
We now present the likelihood-informed subspace dimensionality reduction approach from \cite{spantiniOptimalLowRank}. Due to sparse measurements and the smoothing nature of the forward operator $\mathbf{G}$, the uncertainty reduction relative to the prior typically is largest in a low-dimensional space, also referred to as the likelihood-informed subspace (LIS). The LIS contains directions that maximize the following Rayleigh quotient:
\begin{equation}
    \frac{\mathbf{f}^\top \mathbf{G}^\top\mathbf{\Gamma}_{\text{obs}}^{-1}\mathbf{G}\mathbf{f}}{\mathbf{f}^\top \bm\Gamma^{-1}\mathbf{f}}.\label{eq:RwHat}
\end{equation}
This ratio normalizes an energy defined by the Fisher information matrix, $\mathbf{G}^\top\mathbf{\Gamma}_{\text{obs}}^{-1}\mathbf{G}$, by an energy defined by the prior precision, $\bm\Gamma^{-1}$.  Directions of parameter space which maximize this quotient can be understood as being maximally informed by the data relative to their prior uncertainty. 

Directions maximizing \eqref{eq:RwHat} are given by the leading right and left eigenvectors $\mathbf{v}_i,\mathbf{w}_i\in\mathbb{R}^d$ of the Fisher information and the prior precision
\begin{flalign}    
    \mathbf{G}^\top\mathbf{\Gamma}_{\text{obs}}^{-1}\mathbf{G}\mathbf{v}_i&=\delta_i^2\bm\Gamma^{-1}\mathbf{v}_i,\label{eq:wHat}\\
    \mathbf{w}_i&=\mathbf{\Gamma}^{-1}\mathbf{v}_i,\label{eq:wTilde}
\end{flalign}
where $\mathbf{v}_i$ are normalized by the inner product induced by the prior precision $\mathbf{v}_i^\top\mathbf{\Gamma}^{-1}\mathbf{v}_i=1$ and $\delta_i^2$ are in descending order. 
If square root factorizations of the prior and noise covariances are available, i.e., $\bm\Gamma=\mathbf{S}\mathbf{S}^\top$ and $\bm\Gamma_{\text{obs}}=\mathbf{S}_{\text{obs}}\mathbf{S}_{\text{obs}}^\top$, the generalized eigenvectors $\mathbf{v}_i$ and $\mathbf{w}_i$ may be defined by $\mathbf{v}_i=\mathbf{S}\bm\nu_i$ and $\mathbf{w}_i=\mathbf{G}^\top \mathbf{S}_{\text{obs}}^{-\top}\bm\omega_i\frac{1}{\delta_i}$, where $\bm{\omega}_i$ and $\bm{\nu}_i$ are the left and right singular vectors and $\delta_i$ the corresponding singular values of 
\begin{equation}
    \mathbf{S}_{\text{obs}}^{-1}\mathbf{G}\mathbf{S}=\sum_{i\geq1}\delta_i\bm\omega_i\bm\nu_i^\top.\label{eq:square-root}
\end{equation}
This computational procedure is similar to the `square-root balancing' procedure used in the balanced truncation method for projection-based model reduction and is better conditioned than computing the generalized eigenvalues directly~\cite{Antoulas2005}.

Assuming that both $\mathbf{\Gamma}_{\text{obs}}$ and $\mathbf{\Gamma}$ are full rank, the number of non-zero eigenvalues $\delta_i^2$ is at most $\text{rank}(\mathbf{G})$. In many practical settings, the number of measurements $m$ is low relative to the dimension of the unknown parameter $\mathbf{f}$, leading to a maximum of $m$ non-zero eigenvalues of~\eqref{eq:wHat} and a maximum LIS dimension of $m$. Additionally, when the forward operator $\mathbf{G}$ exhibits smoothing properties, the non-zero eigenvalues of~\eqref{eq:wHat} often exhibit strong decay, enabling a limited number of leading eigenvectors to represent the parameter space where the data are most informative.

We collect the first $r$ eigenvectors corresponding to the largest eigenvalues defined by \eqref{eq:wHat} and \eqref{eq:wTilde} as columns of two matrices 
\begin{flalign}
    \approxBasis &=[\mathbf{v}_1,...,\mathbf{v}_r]\in\mathbb{R}^{d\times r}\label{eq:wHatMatrix},\\
    \proBasis &= [\mathbf{w}_1,...,\mathbf{w}_r]\in\mathbb{R}^{d\times r}.\label{eq:wTildeMatrix}
\end{flalign}
Note that $\mathbf{V}^\top\mathbf{W}=\mathbf{I}_r$, so the matrices $\approxBasis$ and $\proBasis$ form an oblique projector given by
\begin{equation}
    \mathbf{P}=\approxBasis\proBasis^\top,\label{eq:projector}
\end{equation}
with $\mathbf{P}^2=\approxBasis\proBasis^\top\approxBasis\proBasis^\top=\approxBasis\proBasis^\top=\mathbf{P}$. The LIS of dimension $r$ is defined to be $\text{Ran}(\mathbf{P})$. We emphasize that $\mathbf{P}$ is an oblique projector: its range is \textit{not} orthogonal to its nullspace.

We can reduce the dimension of the inverse problem by projecting the unknown parameter $\mathbf{f}$ onto the LIS using $\mathbf{P}$, yielding the following approximation of~\eqref{eq:linearIP}:
\begin{flalign*}
    \mathbf{y}&=\mathbf{G}\mathbf{P}\mathbf{f}+\bm\epsilon \equiv \GSpa \mathbf{f} + \bm\epsilon,\\
    \textcolor{black}{\GSpa}&\textcolor{black}{=\mathbf{CK}^{-1}\mathbf{VW}^\top},
\end{flalign*}
where we have defined $\GSpa\in\mathbb{R}^{m\times d}$. This approximate measurement model defines  the following posterior mean and covariance approximations:
\begin{flalign}
    \bm\mu_{\text{pos}}^{\text{\tiny OLR}}&=\bm\mu+\bm\Gamma \GSpa^\top(\GSpa\bm\Gamma \GSpa^\top+\bm\Gamma_{\text{obs}})^{-1}(\mathbf{y}-\GSpa\bm\mu),\label{eq:posMeanSp}\\
    \bm\Gamma_{\text{pos}}^{\text{\tiny OLR}}&=\bm\Gamma-\bm\Gamma \GSpa^\top(\GSpa\bm\Gamma \GSpa^\top+\bm\Gamma_{\text{obs}})^{-1}\GSpa\bm\Gamma,\label{eq:posCovSp}
\end{flalign}
where OLR stands for ``optimal low-rank''. The work~\cite{spantiniOptimalLowRank} shows that these approximations are optimal in the following senses: first, when $\bm{\mu}=\boldsymbol{0}$, the posterior mean approximation \eqref{eq:posMeanSp} minimizes the Bayes risk $\mathbb{E}\left[\Vert\bm{\mu}_{\text{pos}}-\bm{\mu}_{\text{pos}}^{\text{\tiny OLR}}\Vert_{\mathbf{\Gamma}_{\text{pos}}^{-1}}\right]$ in the class of mean approximations given by a rank-$r$ linear operator applied to $\mathbf{y}$.
Second, 
the posterior covariance approximation \eqref{eq:posCovSp} is the optimal rank-$r$ negative semi-definite update of the prior covariance, i.e. $\{\mathbf{\widehat{\Gamma}}_{\text{pos}}:=\mathbf{\Gamma}-\mathbf{B}\mathbf{B}^\top, \text{with rank}(\mathbf{B})\leq r \}$, as measured by the F\"orstner distance $d_F$, defined as
\begin{equation}
    d_F(\mathbf{\Gamma}_{\text{pos}},\mathbf{\widehat{\Gamma}}_{\text{pos}})=\sqrt{\text{tr}\left[\text{ln}^2(\mathbf{\Gamma}_{\text{pos}}^{-0.5}\mathbf{\widehat{\Gamma}}_{\text{pos}}\mathbf{\Gamma}_{\text{pos}}^{-0.5}) \right]}=\sqrt{\sum_i\text{ln}^2(\lambda_i)}\label{eq:foerstnerMetric}
\end{equation}
where $\lambda_i$ are the generalized eigenvalues of pencil $(\mathbf{\Gamma}_{\text{pos}},\mathbf{\widehat{\Gamma}}_{\text{pos}})$. The metric in \eqref{eq:foerstnerMetric} compares the shape and orientation of the uncertainty encoded in the posterior covariance matrices and is invariant to rotation and scaling. It treats over- and underapproximations similarly by divergence to $\infty$. Moreover, $d_F(\mathbf{\Gamma}_{\text{pos}},\mathbf{\widehat{\Gamma}}_{\text{pos}})=d_F(\mathbf{\Gamma}_{\text{pos}}^{-1},\mathbf{\widehat{\Gamma}}_{\text{pos}}^{-1})$, resulting in a metric not only for the posterior covariance, but simultaneously for the posterior precision. Note that the optimality of the posterior covariance approximation is independent of the data and of the prior mean. 

Note that the approximate optimal-low-rank forward operator $\GSpa=\mathbf{G}\mathbf{P}$ has dimensions $m\times d$ but rank $r$. This low-rank structure \textcolor{black}{can be exploited for both efficient storage and computation of the posterior mean and covariance approximations. However, we emphasize that this OLR dimension reduction approach does \textit{not} result in a reduced \textit{model} as the approach in \Cref{sec:MOR}. This means that downstream computations requiring solution of the model, e.g., to propagate posterior parameter samples to quantities of interest that depend on the system state, require solution of the high-dimensional linear system. In \Cref{sec:LIS}, we use the subspace identified in this optimal-low-rank dimension reduction approach to develop a new projection-based \textit{model reduction} technique. A detailed analysis and comparison of the time and storage costs of the proposed method with existing dimension and model reduction methods is provided in \Cref{sec:costs}, and numerical comparisons are provided in \Cref{sec:Examples}. }

\section{A likelihood-informed subspace method for model reduction}\label{sec:LIS}
In this section, we present a novel model reduction method exploiting the low-dimensionality of the LIS for inference of the forcing of linear static systems. We specifically consider the linear Gaussian inverse problem introduced in \Cref{sec:linGauIP}, where the right-hand-side forcing is inferred from linear measurements of the state, and the prior and noise distributions are Gaussian. We first derive a likelihood-informed subspace (LIS) reduced model by a Petrov-Galerkin projection into the LIS and show how the reduced model is used in the context of linear Bayesian inverse problems \textcolor{black}{in \Cref{sec:LIS_MR}. \Cref{sec:costs} compares the cost and storage requirements of the proposed method to the requirements of the full model, OLR and POD.}

\subsection{\textcolor{black}{Likelihood-informed model reduction}}\label{sec:LIS_MR}
Let $\approxBasis\in\mathbb{R}^{d\times r}$ and $\proBasis\in\mathbb{R}^{d\times r}$ contain the dominant left and right eigenvectors of the Fisher information and the prior precision given in \eqref{eq:wHatMatrix} and \eqref{eq:wTildeMatrix}. 
We use the trial and test basis $\approxBasis$ and $\proBasis$ in a Petrov-Galerkin projection described in \eqref{eq:approximationFullState}-\eqref{eq:reducedPDE} to obtain the reduced model and approximate the observational model given in \eqref{eq:linearPDE} by defining an approximate  model output given by
\begin{flalign}
    \mathbf{y}=\mathbf{\widehat{C}}_{\textLISMR}\mathbf{\widehat{u}}_{\textLISMR}+\bm\epsilon,\quad \text{s.t.:} \quad \mathbf{\widehat{K}}_{\textLISMR}\mathbf{\widehat{u}}_{\textLISMR}=\mathbf{\widehat{f}}_{\textLISMR},\label{eq:LISapproximatedIP}
\end{flalign}
where $\mathbf{\widehat{C}}_{\textLISMR}=\mathbf{C}\approxBasis\in\mathbb{R}^{m\times r}$ maps from the reduced state to the measurements, $\mathbf{\widehat{K}}_{\textLISMR}=\proBasis^\top\mathbf{K}\approxBasis\in\mathbb{R}^{r\times r}$ is the reduced model obtained using the proposed method and $\mathbf{\widehat{f}}_{\textLISMR}=\proBasis^\top\mathbf{f}\in\mathbb{R}^r$ is the projected parameter. Similar to \eqref{eq:PODIP}, we define a reduced inverse problem 
\begin{flalign*}
    \mathbf{y}&=\GhatLIS\mathbf{\widehat{f}}_{\textLISMR}+\bm\epsilon,\\
    \textcolor{black}{\GhatLIS}&\textcolor{black}{=\mathbf{C}\approxBasis(\proBasis^\top\mathbf{K}\approxBasis)^{-1}}
\end{flalign*}
where $\GhatLIS=\mathbf{\widehat{C}}_{\textLISMR}\mathbf{\widehat{K}}_{\textLISMR}^{-1}\in\mathbb{R}^{m\times r}$ is the reduced forward operator and LIS stands for ``likelihood-informed subspace''. The reduced parameter $\mathbf{\widehat{f}}_{\textLISMR}$ has Gaussian prior distribution with $\bm{\widehat{\mu}}_{\textLISMR}=\proBasis^\top\bm\mu\in\mathbb{R}^{r}$ and $\mathbf{\widehat{\Gamma}}_{\textLISMR}=\proBasis^\top\mathbf{\Gamma}\proBasis\in\mathbb{R}^{r\times r}$ leading to a Gaussian posterior defined by
\begin{flalign}
    \bm{\widehat{\mu}}_{\text{pos}}^{\textLISMR}&=\bm{\widehat{\mu}}_{\textLISMR}+\mathbf{\widehat{\Gamma}}_{\textLISMR} \GhatLIS^\top(\GhatLIS\mathbf{\widehat{\Gamma}}_{\textLISMR} \GhatLIS^\top+\mathbf{\Gamma}_{\text{obs}})^{-1}(\mathbf{y}-\GhatLIS\bm{\widehat{\mu}}_{\textLISMR}),\label{eq:posMeanLISReduced}\\
    \mathbf{\widehat{\Gamma}}_{\text{pos}}^{\textLISMR}&=\mathbf{\widehat{\Gamma}}_{\textLISMR}-\mathbf{\widehat{\Gamma}}_{\textLISMR} \GhatLIS^\top(\GhatLIS\mathbf{\widehat{\Gamma}}_{\textLISMR} \GhatLIS^\top+\bm\Gamma_{\text{obs}})^{-1}\GhatLIS\mathbf{\widehat{\Gamma}}_{\textLISMR},\label{eq:posCovLISReduced}
\end{flalign}
where $\bm{\widehat{\mu}}_{\text{pos}}^\textLISMR\in\mathbb{R}^r$ and $\mathbf{\widehat{\Gamma}}_{\text{pos}}^\textLISMR\in\mathbb{R}^{r\times r}$ lie in the low-dimensional subspace. The corresponding high-dimensional approximation of the posterior mean of $\mathbf{f}$ is given by 
\begin{equation}
    \bm\mu_{\text{pos}}^\textLISMR=\textcolor{black}{(\mathbf{I}_d-\approxBasis\proBasis^\top)\bm\mu+}\approxBasis\bm{\widehat{\mu}}_{\text{pos}}^\textLISMR\in\mathbb{R}^d,
\end{equation}
and the low-rank update posterior covariance approximation is given by
\begin{equation}
    \mathbf{\Gamma}_{\text{pos}}^\textLISMR=\mathbf{\Gamma}-\mathbf{\Gamma}\proBasis\mathbf{\widehat{G}}_{\textLISMR}^\top(\mathbf{\widehat{G}}_{\textLISMR}\proBasis^\top\mathbf{\Gamma}\proBasis\mathbf{\widehat{G}}_{\textLISMR}^\top+\mathbf{\Gamma}_{\text{obs}})^{-1}\mathbf{\widehat{G}}_{\textLISMR}\proBasis^\top\mathbf{\Gamma}.\label{eq:posCovLISMapped}
\end{equation}
After calculating $\approxBasis$ and $\proBasis$ once and setting up \eqref{eq:LISapproximatedIP} in the offline phase, we can solve a $r$-dimensional inverse problem using \eqref{eq:posMeanLISReduced}--\eqref{eq:posCovLISMapped} in the online phase. Note that the reduced forward operator $\GhatLIS$ has storage and evaluation cost \textcolor{black}{in} $\mathcal{O}(mr)$, enabling online costs to scale independently of the original dimension $d$. This can lead to significant computational savings for $r\ll d$. \textcolor{black}{If additional downstream computations rely on posterior state samples, we can directly sample the reduced state by solving the reduced PDE instead of using the full model.}

In many engineering applications, we do not have prior information about the unknown parameters in the directions of all degrees of freedom. This results in the singularity of prior covariances and the prior precision is not defined. For these cases, the directions $\approxBasis$ and $\proBasis$ can be interpreted as the eigenvectors solving the nonsymmetric eigenvalue problem $\mathbf{\Gamma}\mathbf{G}^\top\mathbf{\Gamma}_{\text{obs}}^{-1}\mathbf{G}\mathbf{v}_i=\delta_i^2\mathbf{v}_i$ and $\mathbf{G}^\top\mathbf{\Gamma}_{\text{obs}}^{-1}\mathbf{G}\mathbf{\Gamma}\mathbf{w}_i=\delta_i^2\mathbf{w}_i$ respectively. Note that square-root factors $\mathbf{\Gamma}=\mathbf{S}\mathbf{S}^\top$ and $\mathbf{\Gamma}_{\text{obs}}=\mathbf{S}_{\text{obs}}\mathbf{S}_{\text{obs}}^\top$ are still available and $\approxBasis$ and $\proBasis$ can be computed using square-root balancing described in \Cref{sec:Spantini}.

\begin{figure}[H]
    \centering
    \includegraphics[width=\textwidth]{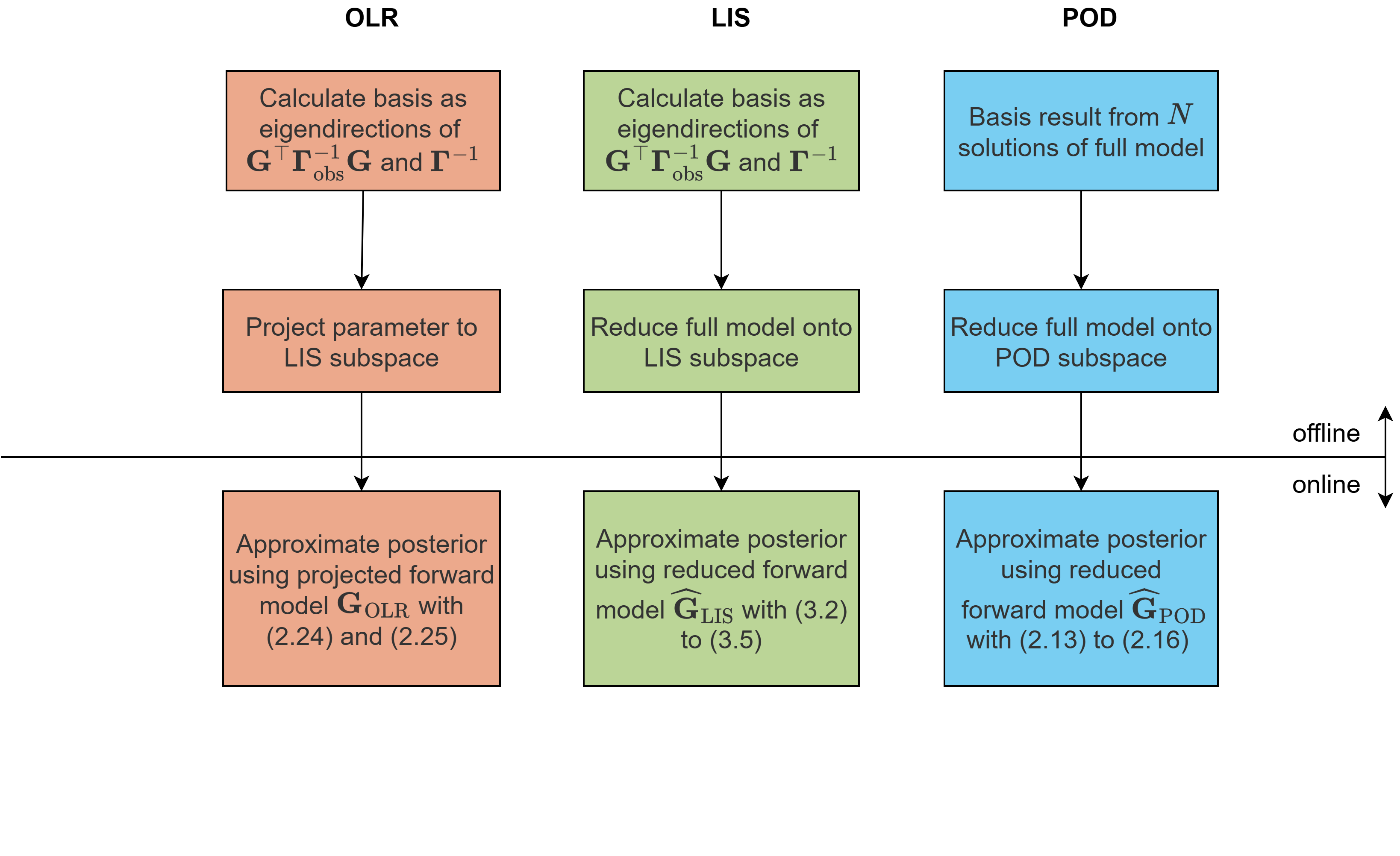}
    \caption{Flowchart diagram of OLR approximation, LIS model reduction and POD model reduction}\label{fig:flowchart}
\end{figure}

 \Cref{algo1} summarizes our likelihood-informed subspace model reduction approach for Bayesian inference of the forcing of linear static systems. \Cref{fig:flowchart} depicts a flowchart containing the key differences between LIS model reduction compared to the OLR approximation and POD model reduction. The offline phase consists of two parts - the basis calculation and the setup of the projected problem. The LIS basis is obtained solving for the eigenvectors of the Fisher information matrix $\mathbf{G}^\top\mathbf{\Gamma}_{\text{obs}}^{-1}\mathbf{G}$ and the prior precision and is therefore entirely based on the inverse problem, whereas the POD basis requires snapshots of the state. To setup the reduced system, OLR approximation only projects the parameter using $\mathbf{P}$ and evaluates the full computational model, while LIS model reduction and POD model reduction define a reduced inverse problem leveraging a reduced-order model. In the online phase, OLR solves $d$-dimensional inverse problems with a projected parameter. LIS model reduction and POD model reduction leverage the reduced inverse problem defined in the offline phase to infer the parameter in $r$ dimensions.

\begin{algorithm}[H]
\caption{Posterior approximation using likelihood-informed model reduction}\label{algo1}
\begin{algorithmic}[1]
\State \textbf{Given:} $\bm\mu,\bm\Gamma,\mathbf{\Gamma}_{\text{obs}},\mathbf{C},\mathbf{K}$
\State Calculate $\approxBasis,\proBasis$ using square-root balancing given in \eqref{eq:square-root}
\State Calculate reduced operator $\mathbf{\widehat{K}}_{\textLISMR}=\proBasis^\top\mathbf{K}\approxBasis$ and $\mathbf{\widehat{C}}_{\textLISMR}=\mathbf{C}\approxBasis$
\State Approximated inverse problem $\GhatLIS=\mathbf{\widehat{C}}_{\textLISMR}\mathbf{\widehat{K}}_{\textLISMR}^{-1}$ with $\mathbf{\widehat{f}}_{\textLISMR}\sim\mathcal{N}(\proBasis^\top\bm\mu,\proBasis^\top\mathbf{\Gamma}\proBasis)$
\State Solve for reduced posterior approximations $\bm{\widehat{\mu}}_{\text{pos}}^\textLISMR$ and $\mathbf{\widehat{\Gamma}}_{\text{pos}}^\textLISMR$
\end{algorithmic}
\end{algorithm}

\subsection{\textcolor{black}{Cost analysis and storage requirements}}\label{sec:costs}
This section provides analysis of the time and storage costs of our method, and compares these to the costs of the POD model reduction and optimal low-rank dimension reduction approaches introduced in \Cref{sec:MOR,sec:Spantini}, as well as to the cost of using the full model assuming direct solvers are available. In our analysis, we consider not only the solution of the linear Bayesian inverse problem through computing the posterior mean and covariance that characterize the uncertain parameter, but also the cost of potential downstream computations including generating system state samples associated with the posterior.
The costs are categorized into a data-independent offline phase including basis computation and approximation of the posterior covariance, a data-dependent online phase in which the posterior mean is approximated, and a downstream computation phase comprising the generation of the state samples. We provide cost analyses for both the case of a dense $\mathbf{K}$ and a sparse banded $\mathbf{K}$, because the latter often arises in structural engineering problems. In what follows, we assume $r\ll d$ and $m\ll d$, and that square-root factors of the prior covariance are available (see Section 3.2.1 of \cite{josieLowRankPrior}). If they have to be computed first, these costs are shared among all approaches and excluded in the comparison. \Cref{tab:dense2} summarizes our analysis of the time costs of the methods. 

If $\mathbf{K}$ is dense, the offline costs of the full model are driven by computing its factorization and are in $\mathcal{O}(d^3)$. If $\mathbf{K}$ has a sparse banded structure, the driving costs are in $\mathcal{O}(d^2m)$ due to the multiplication of $\mathbf{\Gamma}\mathbf{G}^\top$. The online costs are in $\mathcal{O}(d^2)$ independent of the structure of $\mathbf{K}$ as we multiply with $\mathbf{\Gamma}_{\text{pos}}\in\mathbb{R}^{d\times d}$. Posterior state sampling is driven by solving the PDE for parameter samples and is in $\mathcal{O}(d^2)$ for a dense $\mathbf{K}$. If $\mathbf{K}$ has sparse banded structure, the costs drop to be in $\mathcal{O}(d)$. 

When using OLR approximation, the offline costs are in $\mathcal{O}(d^3)$ due to the factorization of the dense $\mathbf{K}$ and solving the singular value decomposition. If $\mathbf{K}$ has sparse banded structure, the costs are in $\mathcal{O}(d^2m)$. Online evaluations can be done in $\mathcal{O}(dr)$ due to multiplication with $\approxBasis$. Equivalently to the full model, posterior state sampling is driven by solving the PDE and is in $\mathcal{O}(d^2)$ for a dense $\mathbf{K}$ and in $\mathcal{O}(d)$ for sparse banded $\mathbf{K}$. 


For the LIS approximation, the offline costs are identical to the ones in OLR and are in $\mathcal{O}(d^3)$ or $\mathcal{O}(d^2m)$, depending on the structure of $\mathbf{K}$. Online costs and downstream costs are dominated by the reconstruction of the high-dimensional posterior mean, covariance, and state samples via multiplication of the reduced mean and covariance with the basis matrix $\approxBasis$, and are thus $\mathcal{O}(dr)$. However, if it is sufficient to compute the reduced quantities (e.g., for downstream sampling only in the reduced representation of the state), then the online and downstream costs would be $\mathcal{O}(r^2)$, leading to greater savings from our approach.


The offline costs of POD approximations are in $\mathcal{O}(d^3+N^2d)$ due to factorization of $\mathbf{K}$ and performing the singular value decomposition of the snapshot matrix. These costs reduce to be in $\mathcal{O}(N^2d)$ if $\mathbf{K}$ has sparse banded structure. The online and downstream costs of POD are identical to LIS and are in $\mathcal{O}(dr)$. If sampling the reduced state is sufficient, for instance when further downstream computations depend only linearly on the state, the costs drop to be in $\mathcal{O}(r^2)$ equivalent to LIS.

Note that LIS and POD are cheaper than OLR downstream computations reducing the costs from $d^2$ to $dr$ for dense $\mathbf{K}$. If sampling the reduced state is sufficient, LIS and POD costs are independent while OLR costs are quadratic or linear in $d$ depending on $\mathbf{K}$.
Moreover, the LIS basis is cheaper compared to POD, if many snapshots (\mbox{$N>\sqrt{dm}$}) are required to obtain accurate POD results ($\mathcal{O}(d^2m)\ll\mathcal{O}(N^2d)$).
\begin{table}[h]
    \centering
    \begin{NiceTabular}{c c|c|c|c|c}
         \CodeBefore
            \rectanglecolor{lightgray}{5-1}{7-6}
         \Body
         & & Full model & OLR & LIS & POD \\
         \toprule
         \Block{3-1}{$\vcenter{\hbox{\rotatebox{90}{dense $\mathbf K$}}}$}
         & offline & $\mathcal{O}(d^3)$ &  $\mathcal{O}(d^3)$ & $\mathcal{O}(d^3)$ & $\mathcal{O}(d^3+N^2d)$\\
         \cmidrule{2-6}
         & online & $\mathcal{O}(d^2)$ & $\mathcal{O}(dr)$ & $\mathcal{O}(dr)$ & $\mathcal{O}(dr)$\\
         \cmidrule{2-6}
         & downstream & $\mathcal{O}(d^2)$ & $\mathcal{O}(d^2)$ & $\mathcal{O}(dr)$ &$\mathcal{O}(dr)$\\
         \midrule
         \Block{3-1}{$\vcenter{\hbox{\rotatebox{90}{banded $\mathbf K$}}}$}
           &  offline &  $\mathcal{O}(d^2m)$ & $\mathcal{O}(d^2m)$ &  $\mathcal{O}(d^2m)$ & $\mathcal{O}(N^2d)$ \\
         \cmidrule{2-6}
         & online & $\mathcal{O}(d^2)$ &  $\mathcal{O}(dr)$ & $\mathcal{O}(dr)$ & $\mathcal{O}(dr)$\\
         \cmidrule{2-6}
         & downstream &$\mathcal{O}(d)$ & $\mathcal{O}(d)$ &  $\mathcal{O}(dr)$ & $\mathcal{O}(dr)$\\
         \bottomrule
    \end{NiceTabular}
    \caption{Total offline (data-independent), online (data-dependent) and downstream costs using full model, OLR, LIS and POD for a dense $\mathbf{K}$ and banded $\mathbf{K}$, where $d$ is the dimension of the parameter, $m$ of the data, $r$ of the reduced model and $N$ the number of snapshots used for the POD basis}
    \label{tab:dense2}
\end{table}
\Cref{tab:storage} lists the key quantities to store for every approach. The full model requires storage of the forward operator $\mathbf{C}$ and $\mathbf{K}$. OLR needs to store the basis $\proBasis$ and $\approxBasis$ on top of $\mathbf{C}$ and $\mathbf{K}$. Note that if posterior state samples are not required, the OLR approximation of the forward operator can be efficiently stored without explicit storage of $\mathbf{K}$. LIS requires storage of basis $\proBasis$ and $\approxBasis$ as well as reduced operators $\mathbf{\widehat{C}}$ and $\mathbf{\widehat{K}}$. As POD is obtained using a Galerkin projection, we only need to store $\bm\Phi_r$ and reduced operators $\mathbf{\widehat{C}}$ and $\mathbf{\widehat{K}}$.
\begin{table}[h]
    \centering
    \begin{NiceTabular}{c c|c|c|c}
          & Full model & OLR & LIS & POD \\
         \toprule
         \Block{2-1}{$\vcenter{\hbox{\rotatebox{0}{basis}}}$}
         & - & $\proBasis\in\mathbb{R}^{d\times r}$ &  $\proBasis\in\mathbb{R}^{d\times r}$ & $\bm\Phi_r\in\mathbb{R}^{d\times r}$\\
         \cmidrule{2-5}
         & - & $\approxBasis\in\mathbb{R}^{d\times r}$ & $\approxBasis\in\mathbb{R}^{d\times r}$ & -\\
         \midrule
         \Block{2-1}{$\vcenter{\hbox{\rotatebox{0}{forward operator}}}$}
          & $\mathbf{C}\in\mathbb{R}^{m\times d}$ & $\mathbf{C}\in\mathbb{R}^{m\times d}$ & $\mathbf{\widehat C}\in\mathbb{R}^{m\times r}$ & $\mathbf{\widehat C}\in\mathbb{R}^{m\times r}$\\
          \cmidrule{2-5}
          & $\mathbf{K}\in\mathbb{R}^{d\times d}$ & $\mathbf{K}\in\mathbb{R}^{d\times d}$ & $\mathbf{\widehat{K}}\in\mathbb{R}^{r\times r}$ &$\mathbf{\widehat{K}}\in\mathbb{R}^{r\times r}$\\
          \bottomrule
    \end{NiceTabular}
    \caption{Storage requirements of full model, OLR, LIS and POD}
    \label{tab:storage}
\end{table}


\section{Examples}\label{sec:Examples}
In this section, we demonstrate the efficacy of our proposed method for two inverse problems drawn from structural mechanics and compare its performance with POD model reduction (\Cref{sec:MOR}) and the OLR approximation presented in \Cref{sec:Spantini}. We apply the developed method to estimate the live load of a cantilever bar in \Cref{sec:bar} and a tunnel with varying ground support in \Cref{sec:tunnel}. The codes used to generate these results are available on GitHub.\footnote{ https://github.com/jakobscheffels/Likelihood-Informed-Model-Reduction-for-Bayesian-Inference-of-Static-Structural-Loads.git}

\subsection{Cantilever bar}\label{sec:bar}
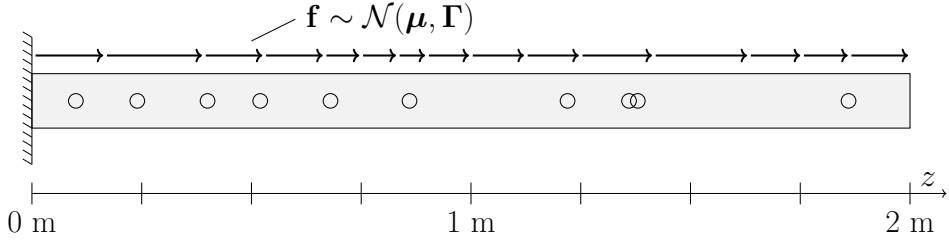
\begin{figure}[H]
    \begin{center}
    \begin{adjustbox}{max width=0.95\linewidth}
    \begin{tikzpicture}

        \draw[thick, fill=gray!10] (0,0) rectangle (24,-1.5);
    
        \draw[thick] (0,1) -- (0,-2.5);
        \foreach \y in {-2.5,-2.25,...,1}
            \draw (0,\y) -- (-0.25,\y+0.15);
    
        
        \draw[thick] (0,-3.3) -- (24,-3.3);
        \foreach \y/\x in {3/0.25,6/0.5,9/0.75,15/1.25,18/1.5,21/1.75}
        \draw[thick] (\y,-3)--(\y,-3.6);
        \foreach \y/\x in {0/0,12/1,24/2} \draw[thick] (\y,-3)--(\y,-3.6) node[below=2pt] {\Huge $\x~\text{m}$};
        \draw[->,thick] (24,-3.3)--(25,-3.3) node[midway,above=2pt]{\Huge $z$};

        \foreach \y in {0.1, 0.24, 0.4, 0.52, 0.68,0.86,1.22,1.36,1.38,1.86}
        \draw[thick] (12*\y,-.75) circle (0.2cm);
    
        \def\z{0.5};
        \draw[->, ultra thick] (0.1,\z) -- (1.95,\z);
        \draw[->, ultra thick] (2.05,\z) -- (4.65,\z);
        \draw[->, ultra thick] (4.75,\z) -- (6.3,\z);
        \draw[->, ultra thick] (6.4,\z) -- (7.95, \z);
        \draw[->, ultra thick] (8.05,\z) -- (8.95,\z);
        \draw[->, ultra thick] (9.05,\z) -- (9.95,\z);
        \draw[->, ultra thick] (10.05,\z) -- (10.75,\z);
        \draw[->, ultra thick] (10.85,\z) -- (11.95,\z);
        \draw[->, ultra thick] (12.05,\z) -- (13.45,\z);
        \draw[->, ultra thick] (13.55,\z) -- (14.95,\z);
        \draw[->, ultra thick] (15.05,\z) -- (16.95,\z);
        \draw[->, ultra thick] (17.05,\z) -- (19.55,\z);
        \draw[->, ultra thick] (19.65,\z) -- (21,\z);
        \draw[->, ultra thick] (21.1,\z) -- (22.3,\z);
        \draw[->, ultra thick] (22.4,\z) -- (23.95,\z);
        
    
        \draw[thick] (6,\z+0.4)--(7.2,1.5);
        \node at (9.8,1.5) {\Huge $\mathbf{f}\sim\mathcal{N}(\bm\mu,\mathbf{\Gamma})$};
        
    \end{tikzpicture}
    \end{adjustbox}
    \end{center}
    \caption{Cantilever bar and circle markings depicting measurement locations}
    \label{fig:Bar}
\end{figure}

Our first example is a cantilever bar of length $L=2~\text{m}$ governed by the static equilibrium equation
\begin{equation*}
    D\frac{d^2u(z)}{dz^2}+q(z)=0,
\end{equation*}
where $z\in\left[0,L\right]$, $D=4\times 10^8 ~\text{N}$ is the constant rigidity, $u(z)$ is the displacement and $q(z)$ the distributed horizontal load modeled as a Gaussian random field. We assume a constant prior mean $\mu_Q =4\times 10^6~\text{N}$ and an exponential auto-correlation given by $\Gamma_{QQ}(z_1,z_2)=\sigma_Q^2\text{exp}\left(-\frac{\vert z_1-z_2\vert}{\theta} \right)$ with standard deviation $\sigma_Q=0.3\cdot \mu_Q$ and correlation length $\theta = L/2$. The boundary conditions are $u(0)=0$, $\frac{du(z)}{dz}|_{z=L}=0$ (see \Cref{fig:Bar}).
The governing equation is solved by the linear finite element method using $n=100$ linear elements.
We discretize the random field using the midpoint method \cite{der1988stochastic} and represent the load within each finite element as $\mathbf{q}\in\mathbb{R}^n$ with $\mathbf{q}_i= q(\Bar{z}_i)$, where $\Bar{z}_i$ is the midpoint of the respective element. The discretized approximation of the random load has a constant mean vector $\bm\mu_Q \in \mathbb{R}^{100}$ where all entries are equal to $\mu_Q$. 
The covariance matrix $\mathbf{\Gamma}_{QQ} \in \mathbb{R}^{100 \times 100}$ is calculated by $\bm\Gamma_{QQ,ij}=\Gamma_{QQ}(\Bar{z}_i,\Bar{z}_j)$ for all elements $i,j$. The nodal force vector $\mathbf{f}\in\mathbb{R}^d$ obtained through integration of the discretized random load over the linear finite element shape functions is given as $\mathbf{f}=\mathbf{L}\mathbf{q}$, where $\mathbf{L}\in\mathbb{R}^{d\times n}$ is the map containing the integration of the linear element shape functions assembled such that the contribution of element $e$ to the $i$-th component of $\mathbf{f}$ is $L_{i,e} = \int_{0}^{L/n} N_i(\chi)\, d \chi$. The shape functions of a linear bar element are given by $N_1(\chi)=1-\frac{\chi\cdot n}{L}$ and $N_2(\chi)=\frac{\chi\cdot n}{L}$ with elemental coordinate $\chi\in[0,\frac{L}{n}]$. The unknown nodal force vector $\mathbf{f}$ therefore is distributed according to a Gaussian prior $\mathcal{N}(\bm\mu,\mathbf{\Gamma})$ with $\bm\mu=\mathbf{L}\bm\mu_Q\in\mathbb{R}^d$ and $\mathbf{\Gamma}=\mathbf{L}\mathbf{\Gamma}_{QQ}\mathbf{L}^\top\in\mathbb{R}^{d\times d}$. 

We consider $m=10$ noisy observations of the displacement degrees of freedom chosen randomly from a uniform distribution along the depth of the bar, with an observation noise of \textcolor{black}{$\mathbf{\Gamma}_{\text{obs}}=\sigma_{\text{obs}}^2\cdot\mathbf{I}_m$}, 
where $\mathbf{I}_m$ is the $m$-dimensional identity matrix and 
\textcolor{black}{$\sigma_{\text{obs}}=0.001$ is set to $5\%$ of the maximum mean displacement. Note that the magnitude of $\sigma_{\text{obs}}$ only influences the singular values $\delta_i$ of $\mathbf{S}_{\text{obs}}^{-1}\mathbf{CK}^{-1}\mathbf{S}$ and not the directions of the likelihood-informed subspace. Larger measurement noise results only in a stronger dependence of the posterior on the prior compared to the measurements; this will decrease the influence of the approximation of the forward model on the overall mean and covariance approximation errors.} 

The location of the observed degrees of freedom are fixed once drawn and are used for all measurement samples (see \Cref{fig:Bar}). 
We generate the measurements through sampling from the prior distribution of $\mathbf{f}$, solving for $\mathbf{u}$, and adding a realization of the observation noise drawn from the distribution $\mathcal{N}(\bm0,\mathbf{\Gamma}_{\text{obs}})$ to the observations $\mathbf{Cu}$. We compute the posterior mean and covariance approximation using LIS model reduction, POD model reduction and OLR dimensionality reduction. As mentioned in \Cref{sec:Spantini}, since we have $m=10$ sparse measurements, we have $\text{rank}(\mathbf{G}^\top\mathbf{\Gamma}_{\text{obs}}^{-1}\mathbf{G})= 10$ leading to $10$ nonzero singular values $\delta^2_i$. The maximum dimension of the LIS is $r=m=10$, and we therefore consider approximations for $r\leq 10$.
 
Because the posterior mean changes depending on the realization of the data, we generate $N_{{\text{rep}}}=200$ realizations of the data $\mathbf{y}^{(i)}$ and report the average of the relative $\ell^2$-norm error between the approximation and the analytical posterior mean obtained with the full order model, given by
 \begin{equation}
    \frac{1}{N_{\text{rep}}}\sum_{i=1}^{N_{\text{rep}}}\frac{\Vert\bm\mu_{\text{pos}}^*(\mathbf{y}^{(i)})-\bm\mu_{\text{pos}}(\mathbf{y}^{(i)})\Vert_2}{\Vert\bm\mu_{\text{pos}}(\mathbf{y}^{(i)})\Vert_2},\label{eq:meanMetric}
\end{equation} 
where $^*$ is OLR, LIS or POD. The posterior covariance approximation is assessed by the F\"orstner distance \eqref{eq:foerstnerMetric} and does not depend on the data. 

\begin{figure}[h]
    \centering
    \includegraphics[width=\linewidth]{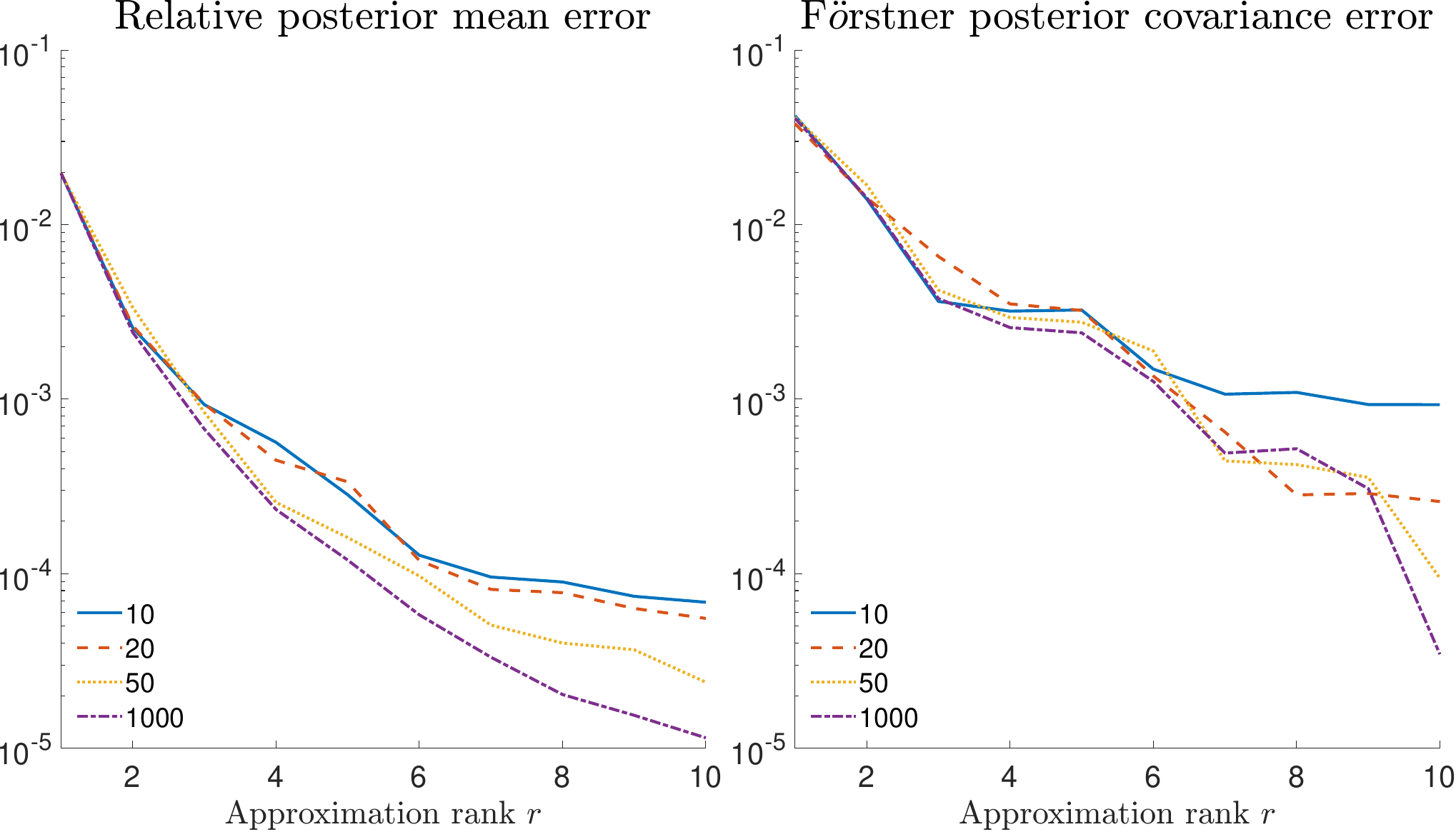}
    \caption{Posterior approximation quality for different number of samples considered in the POD basis calculation for the cantilever bar example}
    \label{fig:errorBarPOD}
\end{figure}

For the POD-reduced model, we perform a study of the number of snapshots needed to compute a POD basis that allows accurate posterior approximation. We construct the POD basis based on solutions of the discretized PDE for samples of $\mathbf{f}$ drawn from the prior distribution. \Cref{fig:errorBarPOD} shows the errors \eqref{eq:meanMetric}, \eqref{eq:foerstnerMetric} of the obtained approximations for increasing number of samples $N$ used to calculate the POD basis from $N=10$ to $N=1000$. Increasing the number of snapshots captures a larger portion of the true variability of the parameter, but converges around a relative posterior mean error of $\approx 10^{-5}$ and covariance error of $\approx 10^{-4}$ for $r=10$. We set $N=10$ for the comparison with LIS model reduction and OLR dimensionality reduction.
 
\begin{figure}[h]
    \centering
    \includegraphics[width=\linewidth]{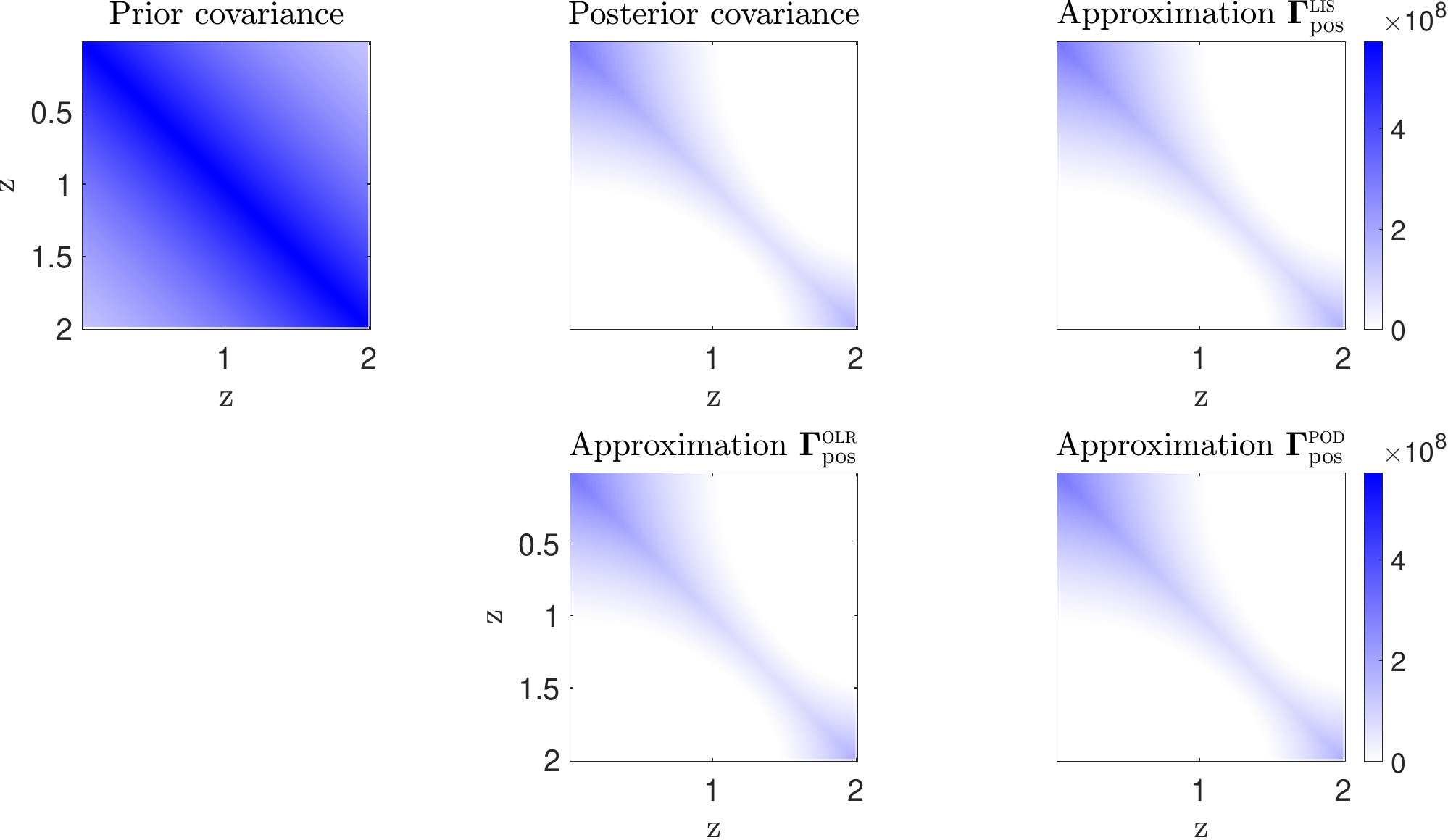}
    \caption{Prior covariance (left), analytical posterior covariance (top center) and the approximation of the posterior covariance for the cantilever bar example obtained using LIS model reduction (top right), OLR dimensionality reduction (bottom center) and POD model reduction (bottom right) for $r=10$}
    \label{fig:covarianceBar}
\end{figure}

 \Cref{fig:covarianceBar} depicts the prior covariance (left), the posterior covariance calculated using the full model (top center) as well as the posterior covariance approximation using LIS model order reduction (top right), OLR dimensionality reduction (bottom center) and POD model reduction (bottom right) for $r=10$. 
We can conclude two main points from \Cref{fig:covarianceBar}. The first point is that uncertainty is reduced in the posterior relative to the prior as expected. The second takeaway is that the solution of the full model and the low-rank approximations are visually identical.  
\begin{figure}[h]
    \centering
    \includegraphics[width=\linewidth]{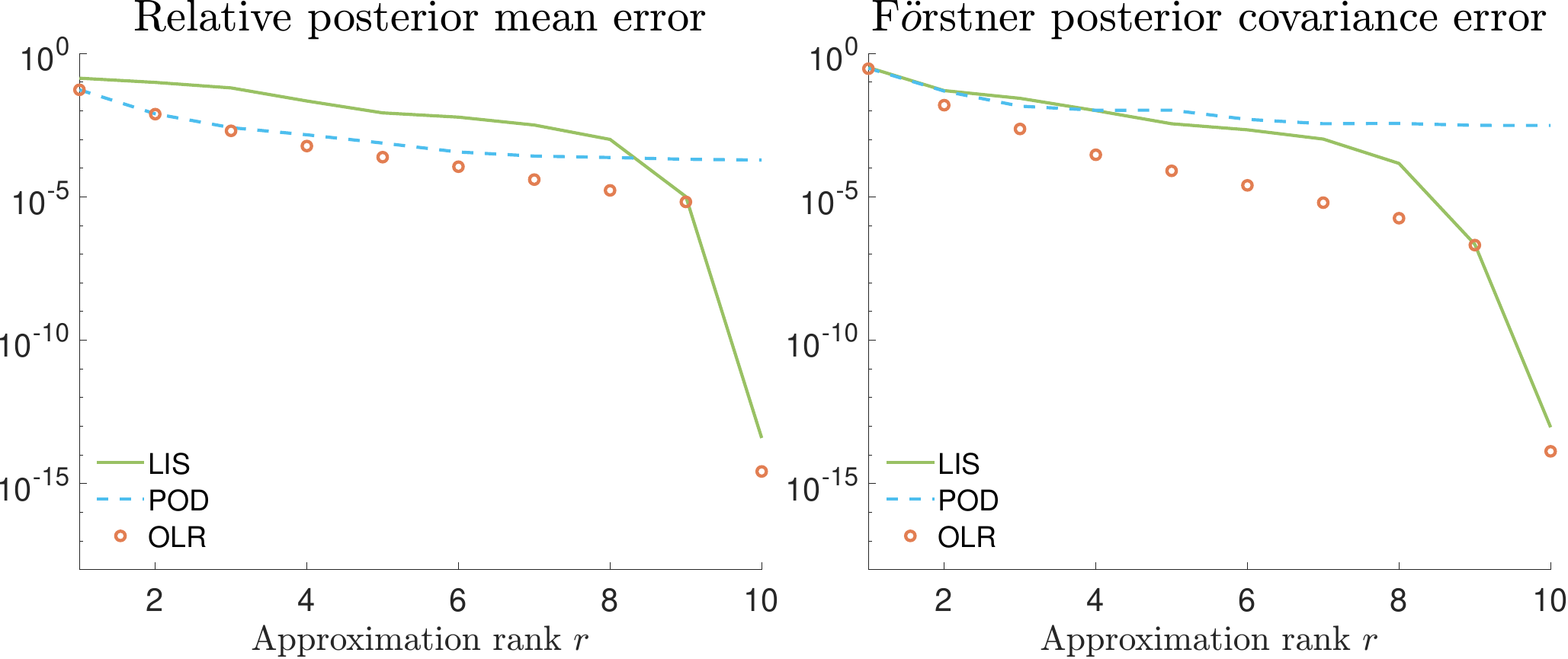}
    \caption{Quantitative assessment of posterior mean (left) and covariance approximation (right) of LIS model reduction, POD model reduction and OLR dimensionality reduction for the cantilever bar example}
    \label{fig:errorMeasuresBar}
\end{figure}

\Cref{fig:errorMeasuresBar} shows the decrease of the posterior mean and covariance error of the three approximations with increasing $r$. For $r=10$, LIS model reduction results in a posterior mean error nine orders of magnitude lower than that of POD, and a posterior covariance error five orders of magnitude lower than that of POD. \textcolor{black}{As illustrated in \Cref{fig:errorBarPOD}, increasing the number of snapshots $N$ to compute the POD basis does not lead to accuracy comparable to LIS.} The LIS model reduction approach errors are near those of the optimal low-rank approximation. We can see that LIS model reduction and OLR dimensionality reduction recognize the intrinsic low-dimensionality of the update as $m=10$ and the respective error drops to numerical precision for $r=10$, while POD model reduction fails to capture this property. As expected, the optimal low-rank aproximation approach results in the best approximations, since the approach uses the high-dimensional stiffness matrix to solve the PDE. LIS model reduction leads to computational savings as the approach computes the posterior approximation by solving a $10$-dimensional reduced model, whereas the dimensionality reduction requires the solution of the original $100$-dimensional finite element discretization. \textcolor{black}{Since the additional cost savings for $r<10$ instead of $r=10$ are negligible relative to the reduction already achieved compared to using the full model as well as the error drop in the approximation quality from $r=9$ to $r=10$, there is little practical justification for using $r<10$. If additional downstream computations depending on posterior state samples are of interest, LIS can introduce additional computational savings according to \Cref{sec:costs}.}

\subsection{Tunnel}\label{sec:tunnel}
\begin{figure}[h]
    \begin{center}
    \begin{adjustbox}{max width=0.95\columnwidth}
    \begin{tikzpicture}
        \def\z{3}; 
        \def\t{2}; 
        \def\r{0.6}; 
        \def\u{-6}; 
        \draw[thick, fill=gray!10] (0,0) rectangle (24,-\z);
        \node[fill=white, inner sep=2pt] at (12,-\z/2) {\Huge Tunnel};
        
        \draw[thick] (0,-\z) -- (0,-\z-\t);
        \draw[thick] (12,-\z) -- (12,-\z-\t);

        \draw[pattern=dots, pattern color=gray]
            (0,-\z) rectangle +(12,-\t);
        \draw[thick,color=white] (0,-\z-\t)--(12,-\z-\t);
        
        \node[fill=white, inner sep=2pt] at (2.3,-\z-\t/2) {\Huge Silty Sand};
        
        \foreach \y in {0.25,0.75,...,8.75}
            \pgfmathsetmacro{\randnum}{rnd}
            \draw[->, ultra thick] (\y,1.5+0.5*\randnum) -- (\y,0.2);
    
        \foreach \y in {9.25,9.75,...,14.75}
            \pgfmathsetmacro{\randnum}{rnd}
            \draw[->, ultra thick] (\y,0.75+0.5*\randnum) -- (\y,0.2);
    
        \foreach \y in {15.25,15.75,...,23.75}
            \pgfmathsetmacro{\randnum}{rnd}
            \draw[->, ultra thick] (\y,1.5+0.5*\randnum) -- (\y,0.2);
            
        \draw[thick] (24,-\z) -- (24,-\z-\t);
        \draw[thick] (12,-\z) -- (12,-\z-\t);

        \draw[pattern=north east lines, pattern color=gray]
            (12,-\z) rectangle +(12,-\t);
        \draw[thick,color=white] (12,-\z-\t)--(24,-\z-\t);
        
        \node[fill=white, inner sep=2pt] at (21.5,-\z-\t/2) {\Huge Mucky Clay};
        
        \node[fill=white, inner sep=2pt] at (12,1.4) {\Huge random live load};

        \draw[thick] (0,\u) -- (24,\u);
        \foreach \y/\x in {3/25,6/50,9/75,15/125,18/150,21/175}
        \draw[thick] (\y,\u+0.3)--(\y,\u-0.3);
        \foreach \y/\x in {0/0,12/100,24/200}
        \draw[thick] (\y,\u+0.3)--(\y,\u-0.3) node[below=2pt] {\Huge $\x~\text{m}$};
        \draw[->,thick] (24,\u)--(25,\u) node[midway,above=2pt]{\Huge $z$};

        \foreach \y in {5.75, 23,26.5,48,65.5,108.25,121.75,133.25,146.75,172}
        \draw[thick] (24*\y/200,-\z+0.75) circle (0.2cm);
        
    \end{tikzpicture}
    \end{adjustbox}
    \end{center}
    \caption{Tunnel with varying ground support subjected to random load and circle markings depicting measurement locations}
    \label{fig:tunnel}
\end{figure}
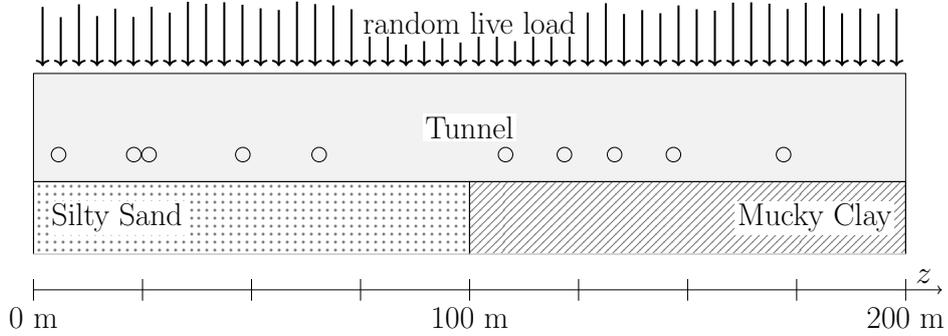

In the second application, we consider a model of a subway tunnel under uncertain loading (\Cref{fig:tunnel}). The tunnel of length $L=200~\text{m}$ is governed by the differential equation
\begin{equation*}
    \frac{d^4u(z)}{dz^4}+\frac{k(z)D u(z)}{\zeta EI}=\frac{q(z)D}{\zeta EI},
\end{equation*}
where $z\in[0,L]$, $u(z)$ is the settlement, $k(z)$ the varying ground support, $D=6.2~\text{m}$ the outer diameter of the hollow beam with thickness $t=0.35~\text{m}$, $\zeta=1/7$ the reduction factor to account for joints between rings, $E=35\times 10^6~\text{kN/m}$ the Young's modulus, $I=\frac{\pi}{64}(D^4-(D-2t)^4)$ the bending moment of inertia and $q(z)$ the live load modeled as a Gaussian random field (see \Cref{fig:tunnel}). We assume $\mu_Q=3~\text{kN}$ and an exponential auto-correlation kernel with $\sigma_Q=\mu_Q$ and $\theta=L/2$. The parameters of the tunnel are taken from \cite{huangTunnel}, where the ground support is modeled with the Winkler model with $k_{\text{Sand}}=33,000~\text{kN/m}^3$ and $k_{\text{Clay}}=5,000~\text{kN/m}^3$. The boundary conditions are $\frac{d^2u(z)}{d^2z}|_{z=0}=0$, $\frac{d^3u(z)}{d^3z}|_{z=0}=0$, $\frac{d^2u(z)}{d^2z}|_{z=L}=0$ and $\frac{d^3u(z)}{d^3z}|_{z=L}=0$.
 \FloatBarrier
We discretize the beam using $n=800$ Bernoulli beam finite elements resulting in $d=1602$ degrees of freedom and use the midpoint method to approximate $\mathbf{q}\in\mathbb{R}^n$ with $\mathbf{q}_i= q(\Bar{z}_i)$ within each element, as described in \Cref{sec:bar}. The discretized approximation of the load has constant mean vector $\bm\mu_Q\in\mathbb{R}^n$ where all entries are equal to $\mu_Q$. The covariance matrix $\mathbf{\Gamma}_{QQ}\in\mathbf{R}^{n\times n}$ is calculated by $\mathbf{\Gamma}_{QQ,ij}=\Gamma_{QQ}(\Bar{z}_i,\Bar{z}_j)$ for all elements $i,j$. The respective nodal force vector $\mathbf{f}\in\mathbb{R}^d$ is given by $\mathbf{f}=\mathbf{L}\mathbf{q}$ and has Gaussian prior with mean $\bm\mu=\mathbf{L}\bm\mu_Q$ and covariance $\mathbf{\Gamma}=\mathbf{L}\mathbf{\Gamma}_{QQ}\mathbf{L}^\top$, where $\mathbf{L}\in\mathbb{R}^{d\times n}$ is the map containing the integration of the cubic element shape functions and is assembled such that the contribution of element $e$ to the $i$-th component of $\mathbf{f}$ is $L_{i,e} = \int_{0}^{L/n} N_i(\chi)\, d \chi$. The four cubic shape functions of the beam element are given by $N_1(\chi)=1-3\left(\frac{\chi\cdot n}{L}\right)^2+2\left(\frac{\chi\cdot n}{L}\right)^3$, $N_2(\chi)=\chi\left[1-2\left(\frac{\chi\cdot n}{L}\right)+\left(\frac{\chi\cdot n}{L}\right)^2\right]$, $N_3(\chi)=3\left(\frac{\chi\cdot n}{L}\right)^2-2\left(\frac{\chi\cdot n}{L}\right)^3$ and $N_4(\chi)=\chi\left[\left(\frac{\chi\cdot n}{L}\right)^2-\frac{\chi\cdot n}{L} \right]$ with elemental coordinate $\chi\in[0,\frac{L}{n}]$, where $N_1$ and $N_3$ describe the elemental settlements and $N_2$ and $N_4$ the elemental rotation at the respective nodes. Note that this discretization naturally defines a singular prior covariance with $\text{rank}(\mathbf{\Gamma})=800<1602=d$. 

Similar to the previous problem (\Cref{sec:bar}), we consider $m=10$ noisy observations of the translational degrees of freedom (settlements) chosen randomly from a uniform distribution along the length of the beam and an observation noise of \textcolor{black}{$\mathbf{\Gamma}_{\text{obs}}=\sigma_{\text{obs}}^2\cdot\mathbf{I}_m$ with $\sigma_{\text{obs}}=0.05$ as $5\%$ of the maximum mean displacement.}  
The locations of the observed degrees of freedom are fixed once drawn and kept for all measurement samples (see \Cref{fig:tunnel}). 
We draw realizations from the distribution of the measurements defined by the prior distribution of $\mathbf{f}$ and the observation noise. 
The quality of the posterior mean approximation is assessed by \eqref{eq:meanMetric} for $N_{\text{rep}}=200$ independent sample observations. 
The posterior covariance approximations is assessed with \eqref{eq:foerstnerMetric} and does not depend on data. 

\begin{figure}
    \centering
    \includegraphics[width=\linewidth]{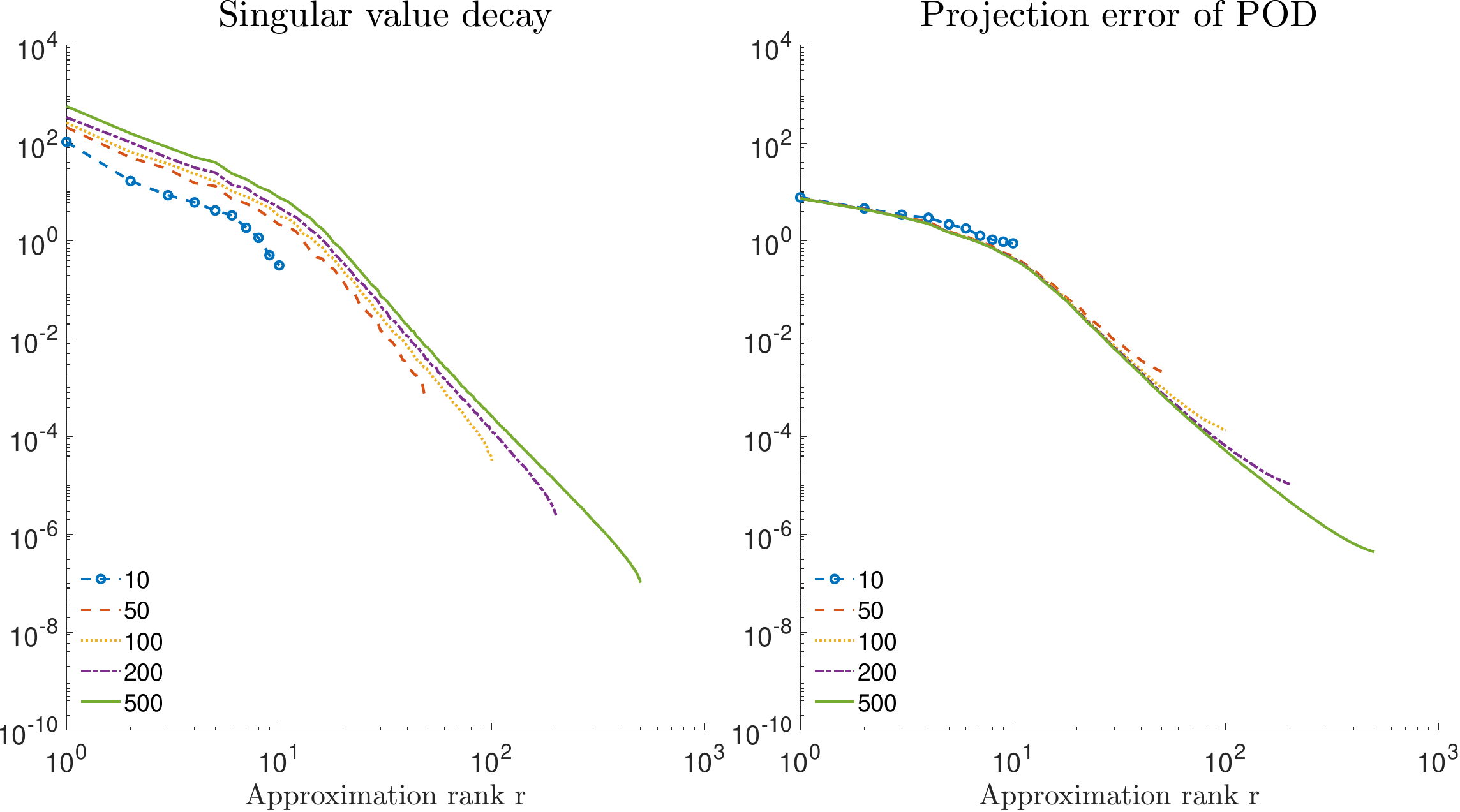}
    \caption{Singular value decay (left) and projection error (right) for different number of samples considered in the POD basis calculation for the tunnel example}
    \label{fig:errorTunnelPOD}
\end{figure}

We again perform a study of the number of snapshots needed to compute a POD basis that yields accurate posterior approximation. Similar to \Cref{sec:bar}, the snapshots are generated by solutions $\mathbf{u}^{(i)}$ of the discretized PDE for samples of $\mathbf{f}^{(i)}$ drawn from its prior distribution for $i=1,...,N$. \Cref{fig:errorTunnelPOD} shows the singular value decay (left) and the projection error $\Vert\mathbf{u}^{(j)}-\bm\Phi_r\bm\Phi_r^\top\mathbf{u}^{(j)}\Vert_2$ (right) of POD-reduced models for increasing number of samples $N$ used to calculate the POD basis from $N=10$ to $N=500$ over a test set of $200$ samples $\mathbf{u}^{(j)}$. The singular values represent the energy of the POD basis and the sum of neglected singular values gives an upper error bound of the projection error of the reduced model. The projection error is a common measure to assess the performance of reduced models by comparing the projected and reconstructed state to the full state. Since $m=10$, we are only interested in the performance for $r\le10$. Based on \Cref{fig:errorTunnelPOD}, we conclude that $N=10$ samples are enough to build the POD-reduced model as more samples do not improve the quality in that range.  

\begin{figure}[h]
    \centering
    \includegraphics[width=\linewidth]{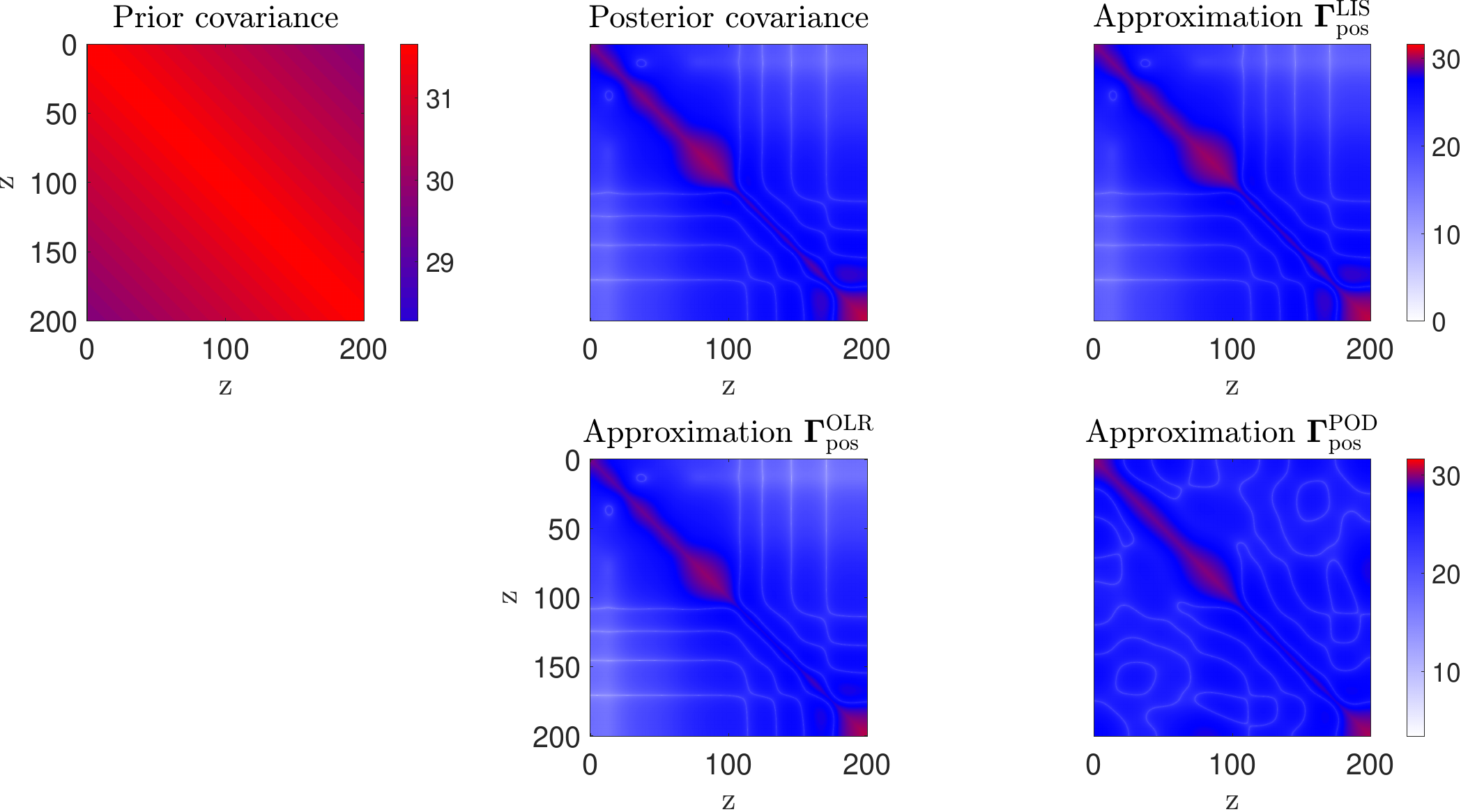}
    \caption{Logarithm of absolute value of prior covariance (left), analytical posterior covariance (top center) and the approximation of the posterior covariance of the translational degrees of freedom of the tunnel example obtained using LIS model reduction (top right), OLR dimensionality reduction (bottom center) and POD model reduction (bottom right) for $r=10$}
    \label{fig:priorTunnel}
\end{figure}

\Cref{fig:priorTunnel} depicts the \textcolor{black}{logarithm of the absolute value of the} prior covariance (left) as well as the posterior covariance of the translational degrees of freedom obtained using the high-dimensional computational model (top center) and the posterior covariance approximation calculated using LIS model reduction (top right), OLR approximation (bottom center) and POD model reduction (bottom right) for $r=10$. We can see that the LIS and OLR approximation visually match the solution of the full model, while there is a difference in the POD approximation \textcolor{black}{failing to accurately capture the uncertainty reduction of the inference}.

\begin{figure}[h]
    \centering
    \includegraphics[width=\linewidth]{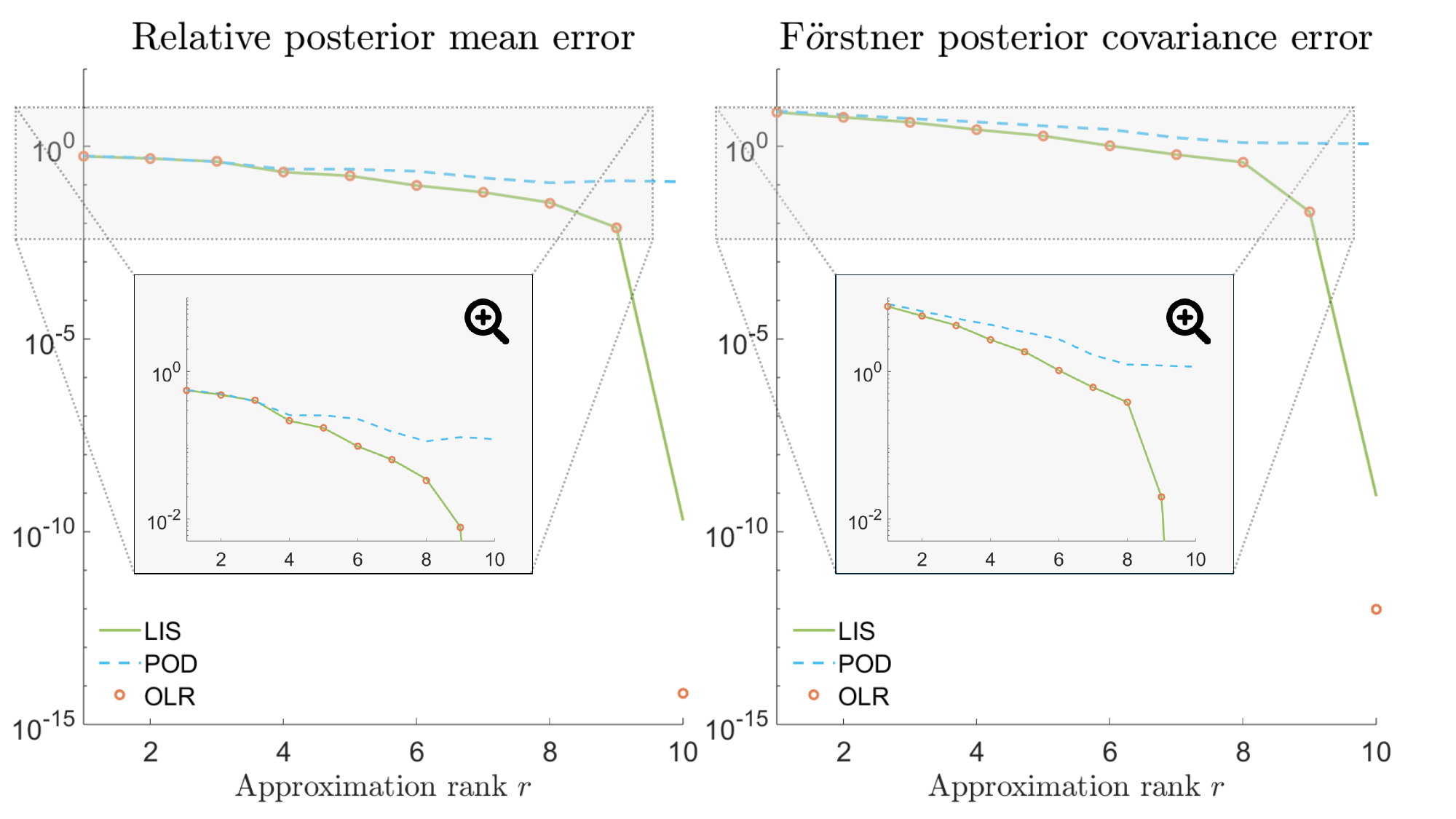}
    \caption{Posterior mean (left) and covariance (right) convergence for increasing $r$ of LIS model reduction, POD model reduction and OLR dimensionality reduction for the tunnel example}
    \label{fig:errorTunnel}
\end{figure}
\Cref{fig:errorTunnel} depicts the decrease of the error of the posterior mean and covariance approximation of LIS, POD and OLR with increasing $r$. The error of LIS approximations matches the quality of the optimal approximation of OLR almost perfectly. Similar to the previous example presented in \Cref{sec:bar}, both methods recognize the dimensionality of the inverse problem as $m=10$ and the error drops significantly for $r=m$. Note that LIS model reduction solves the inverse problem in a $r$-dimensional space, whereas OLR uses the $1602$-dimensional space. For $r=10$, the reduced model needs less than $1\%$ of the degrees of freedom of the full-order model to achieve similar results. 

The error of the POD approximation stays almost constant for $r\le10$, while the error in the LIS approximations decrease multiple orders of magnitude ($\approx 22\times 10^1$ to $2\times10^{-9}$). An explanation is that the modes of the POD basis do not consider the map between the parameter and the sparse observations but try to reconstruct the full state. For accurate posterior approximation, reconstruction of the full state is not necessary, but exploitation of the low-dimensional subspace in the transition from prior to posterior is, which POD fails to capture.

\section{Concluding Remarks}\label{sec:Conclusion}
This paper presents a new likelihood-informed subspace model reduction method for linear systems with right-hand-side uncertainties that enables the efficient solution of \textcolor{black}{linear} Bayesian inverse problems \textcolor{black}{and downstream computations based on posterior state samples}, a key task in model updating and calibration for digital twins. The reduced-order model is obtained by projecting the discretized PDE into the likelihood informed subspace, the space spanned by the dominant generalized eigenvectors of the Fisher information and the prior precision. We apply the developed method to two examples from structural engineering and compare the approximation quality with reduced-order models obtained using POD and OLR approximation. In this comparison, we show that the developed method achieves almost identical results as OLR approximation for reduced dimension equal to the number of measurements, while reducing the inverse problem to the dimension of the update from the prior to posterior distribution. Furthermore, we demonstrate that the tailored reduced model outperforms the generic POD-based model, achieving posterior errors multiple orders of magnitude lower than POD while also requiring fewer computations in the offline phase because snapshot data are not required for our approach.

\textcolor{black}{Future work in this direction will build on the foundation we have established for linear inverse problems to develop likelihood-informed subspace model reduction for nonlinear inverse problems. For example, one could consider extensions of the method} for inference of material parameters resulting in left-hand-side uncertainties and nonlinear Bayesian inverse problems, where no analytical solution for the posterior distribution is available and its properties have to be assessed numerically. A model reduction method for nonlinear inverse problems building on the method presented here might build upon the extension of the ideas of \cite{spantiniOptimalLowRank} to the nonlinear setting presented in \cite{cuiLikelihoodinformedDimensionReduction2014b}. Additionally, one could consider the structural reliability problem as source of nonlinearity to try to identify failure events. In this context, one can use the ideas presented above to develop a model reduction method exploiting the so-called failure informed subspace \cite{uribe2021cross} to accurately approximate failure probabilities using low-rank models. Open questions also remain regarding analysis of the errors of the approximations presented here. 

\subsection*{Declarations}
\textbf{Conflict of interests} \quad The authors declare no conflict of interest. EQ is an Editorial Board member of the special issue on Reduced Order Modeling, Generative AI, and SciML in Digital Twins.

\noindent
\textbf{Data availablility} \quad The code used to generate the results presented above can be downloaded from our GitHub: https://github.com/jakobscheffels/Likelihood-Informed-Model-Reduction-for-Bayesian-Inference-of-Static-Structural-Loads.git

\noindent
\textbf{Replication of results} \quad  The data can be generated using the MATLAB code from our GitHub: https://github.com/jakobscheffels/Likelihood-Informed-Model-Reduction-for-Bayesian-Inference-of-Static-Structural-Loads.git

\noindent
\textbf{Funding} \quad This work was funded by the TUM Georg Nemetschek Institute and the Technical University of Munich - Institute for Advanced Study, Germany. 
EQ also acknowledges support from the United States Air Force Office of Scientific Research under
FA9550-24-1-0105 (program manager Dr.\ Fariba Fahroo).

\noindent
\textbf{Author contributions} \quad JS, EQ, IP and EU contributed to the development of the presented method. Implementation and experiments were run by JS. The first draft of the manuscript was written by JS before being edited by all authors. The final manuscript was read and approved by all authors.

\noindent
\textbf{Ethics approval and consent to participate} \quad Not applicable.




\printbibliography
\end{document}